\documentclass[12pt]{article}
\usepackage{amsmath}
\usepackage{amssymb}

\title{Majorants of meromorphic functions \\ with fixed poles}
\author{A. D. Baranov, A. A. Borichev, V. P. Havin}
\date{}

\begin{document}
\maketitle
\sloppy
\noindent
{\bf Abstract.}
Let $B$ be a meromorphic Blaschke product in the upper half-plane
with zeros $z_n$ and let $K_B=H^2\ominus BH^2$ be the 
associated model subspace of the Hardy class.
In other words, $K_B$ is the space of square summable 
meromorphic functions with the poles at the points $\overline z_n$.
A nonnegative function $w$ on the real line is said to be 
an admissible majorant for $K_B$ if there is a non-zero
function $f\in K_B$ such that $|f|\le w$
a.e. on $\mathbb{R}$.
We study the relations between the 
distribution of the zeros of a Blaschke product $B$ and the class 
of admissible majorants for the space $K_B$. 

\bigskip

\begin{center}
{\bf {\large Introduction}\footnote{A.D. Baranov and V.P. Havin are supported 
by RFBR grants No 03-01-00377 and No 06-01-00313.}}
\end{center}

Let $a>0$ and let $PW_a$ be the Paley--Wiener space of entire functions
of exponential type at most $a$, whose restrictions 
to the real axis $\mathbb{R}$ belong to $L^2(\mathbb{R})$. 
It is well known that the space $PW_a$
coincides with the Fourier image of the space of square integrable 
functions supported in the interval $(-a,a)$.
\smallskip

A nonnegative function $w$ on the real axis $\mathbb{R}$ is said
to be an admissible majorant for the Paley--Wiener space $PW_a$ if
there exists a nonzero function $f\in PW_a$ such that
$|f(x)|\le w(x)$ almost everywhere on $\mathbb{R}$.
It is an important problem of harmonic analysis to describe
the class of admissible majorants for the Paley--Wiener spaces.
An obvious necessary condition is the convergence of the 
logarithmic integral, that is,
\begin{equation}
\label{1}
{\cal L}(w)=\int\limits_\mathbb{R}
\frac{\log^+ w^{-1}(x)}{1+x^2}dx<\infty.
\end{equation}
A sufficient condition of the admissibility is given by the famous
Beurling--Malliavin theorem \cite{bm}: 
{\it if $w$ satisfies (\ref{1}) and
the function
$$
\Omega=-\log w
$$
is Lipschitz on $\mathbb{R}$, then $w$ is an admissible majorant for any space
$PW_a$, $a>0$.} This is a very deep result and quite a few proofs 
are known now (see \cite{hj, ko2, ko3}). 
\smallskip

Let us also mention another, much simpler, sufficient condition:
{\it if $w$ is an even function decreasing on $\mathbb{R}_+ = [0,\infty)$, 
and ${\cal L}(w)<\infty$, then the majorant
$w$ is admissible for any space $PW_a$, $a>0$.}
\smallskip

A new approach to the Beurling--Malliavin theorem was recently
proposed by V.P. Havin and J. Mashreghi \cite{hm1,hm2}.
This approach is based on the study of the Hilbert transform 
of the function $\Omega$. Combined with a recent development of the same 
authors and F. Nazarov \cite{hmn} this approach yields a new (and, probably,
the shortest) proof of the Beurling--Malliavin theorem.
\smallskip

Another advantage of this approach is that it is applicable to a certain
class of spaces of analytic functions, which generalize 
the Paley--Wiener spaces, namely, to the so-called model subspaces 
of the Hardy class. 
\smallskip

Let $\Theta$ be an inner function in the upper half-plane
$\mathbb{C}^+$, that is, a bounded analytic function such that
$\lim\limits_{y\to 0+}|\Theta(x+iy)|=1$ for almost all
$x\in\mathbb{R}$ with respect to Lebesgue measure. 
With an inner function $\Theta$ we associate the model subspace
$$
K_{\Theta}=H^2\ominus \Theta H^2
$$ 
of the Hardy class $H^2$ in the upper half-plane.
These subspaces (and their analogs for the unit disc) 
play an outstanding role both in function and operator theory
(see \cite{cr, nik, nk12}), in particular, in the Sz.-Nagy--Foias
model for contractions in a Hilbert space. It is well known that 
any subspace of $H^2$ coinvariant with respect
to the semigroup of shifts $(U_t)_{t\ge 0}$, $U_tf(x)=e^{itx}f(x)$,
is a ${K_{\Theta}}$ for a certain
inner function $\Theta$.
\smallskip

We mention two important particular cases
of the model subspaces. If $\Theta(z)=\exp(iaz)$, $a>0$, then
${K_{\Theta}}=\exp(iaz/2)PW_{a/2}$. On the other hand, 
if $B$ is a Blaschke product
with zeros $z_n$ of multiplicities $m_n$, 
that is,
$$
B(z)=\prod\limits_n e^{i\alpha_n}
\left(\frac{z-z_n}{z-\overline z_n}\right)^{m_n}
$$
(here $\alpha_n \in \mathbb{R}$ and the factors $e^{i\alpha_n}$ 
ensure the convergence of the product),
then the subspace $K_B$ admits a simple geometrical description:
it coincides with the closed linear span in 
$L^2$ of the fractions $(z-\overline z_n)^{-k}$, $1\le k\le m_n$.
\smallskip

The approach of Havin and Mashreghi makes it possible to obtain
analogs of the Beurling--Malliavin theorem for the model subspaces.
We say that $w$ is {\it an admissible majorant for the space} $K_\Theta$,
if there exists a nonzero function $f\in K_\Theta$ such that
$|f(x)|\le w(x)$ almost everywhere on $\mathbb{R}$.
The class of admissible majorants for $K_\Theta$ we denote by
${\rm Adm}(\Theta)$. Note that the condition ${\cal L}(w)<\infty$
is necessary for the inclusion $w\in {\rm Adm}(\Theta)$ for any
inner function $\Theta$, since
$$
\int\limits_\mathbb{R} \frac{\log |f(x)|}{1+x^2}dx>-\infty
$$
for any nonzero $f\in H^2$ (see, e.g., \cite{hj}, p. 32--36).
\smallskip

In \cite{hm1,hm2} the authors consider mainly the case
where the inner function $\Theta$ is meromorphic in the whole complex plane. 
Then, up to a unimodular constant, $\Theta$ is of the form
$$
\Theta(z)=\exp(iaz)B(z),
$$
where $a\ge 0$ and $B$ is a Blaschke product with zeros tending to infinity.
In this case there is a well-defined branch of the argument 
of $\Theta$ on $\mathbb{R}$, that is, there exists an increasing
function $\varphi$ such that $\Theta(t)=\exp(i\varphi(t))$, $t\in\mathbb{R}$.
Moreover, $\Theta'(t) = i\varphi'(t) \Theta(t)$ and
\begin{equation}
\label{1a}
\varphi'(t)=|\Theta'(t)|=a+2\sum\limits\frac{m_n {\rm Im}\, z_n}{|t-z_n|^2},\quad t\in\mathbb{R}.
\end{equation}
If $B$ is a meromorphic Blaschke product,
then the model subspace $K_B$ may be interpreted as the space of meromorphic
functions with fixed poles in the lower half-plane, which are 
square summable on $\mathbb{R}$.
\smallskip

The sufficient conditions of admissibility obtained in \cite{hm1,hm2} 
are expressed in terms of the Hilbert transform of the function $\Omega$
and the argument $\varphi$ of the inner function $\Theta$. We state
some of these results in \S 2. The structure of the class ${\rm Adm}(\Theta)$
is especially well understood in two model situations. 
In the first case, which is considered in \cite{hm2}, 
it is assumed
that the argument of the inner function grows almost linearly, that is,
\begin{equation}
\label{2}
C_1 \le \varphi'(t)\le C_2, \qquad t\in\mathbb{R},
\end{equation}
for some positive constants $C_1$ and $C_2$.
In this case the class of admissible majorants essentially coincides with
the class of admissible majorants for the Paley--Wiener space
(we give a precise statement in \S 2). On the other hand, if
the zeros of $B$ are in a sense sufficiently sparse near
the real axis (for example, if they are situated on the ray 
$\{z=iy: y>0\}$), then there exists
a "minimal" positive admissible majorant for $K_B$
(see \cite{hm1} and, also, \cite{bar04}).
\smallskip

There is a certain gap between these two cases. The goal of this paper
is to fill in this gap and to show how the class of admissible majorants 
${\rm Adm}(B)$ depends on the distribution of zeros of the Blaschke product
$B$. We study admissible majorants for a class of model subspaces generated 
by meromorphic Blaschke products with regularly distributed zeros.
In particular, we consider in detail
the case when the zeros lie in a strip or in a half-strip and have 
a power growth.
\smallskip

We will frequently work with interpolating Blaschke products $B$
with zeros $z_n$. In this case the functions $f$ in $K_B$ may be 
characterized by the representation
\begin{equation}
\label{00}
f(z)=\sum\limits_n\frac{c_n}{z-\overline z_n},
\end{equation}
where $\sum_n ({\rm Im}\, z_n)^{-1}|c_n|^2<\infty$ (see \cite{nik}), and thus the problem 
reduces to the study of majorants for the series of the form (\ref{00}).
This representation of the elements of $K_B$
makes it possible to relate our problem
to recent results on quasianalyticity \cite{ko01}
and weighted polynomial approximation \cite{bs}.
\smallskip

In what follows we make use of
the following notations:
given nonnegative functions $g$ and $h$ we write
$g\asymp h$ if $C_1 h\le g\le C_2 h$
for some positive constants $C_1$ and $C_2$
and for all admissible values of the variables.
Letters $C$, $C_1$, etc. will denote various constants
which may change their values in different occurrences.
\bigskip

\begin{center}
{\bf \large \S 1. Main results}
\end{center}

As in \cite{hm1,hm2}, in the present paper we restrict ourselves 
to the case of meromorphic inner functions.
Let us consider the following example. 
Let $B$ be the 
Blaschke product with simple zeros $z_n=n+i$, $n\in\mathbb{Z}$.
Then, clearly, the argument of $B$ satisfies (\ref{2}).
Moreover, it is easy to show (see \S 2) that 
${\rm Adm}(B)={\rm Adm}(e^{2\pi iz})$.
Now let us take only a half of the zeros: consider the Blaschke product
$B_1$ with zeros at the points $z_n=n+i$, $n\in\mathbb{N}$. It is a natural question
whether there is a qualitative difference between the classes ${\rm Adm}(B_1)$
and ${\rm Adm}(B)$. Clearly, ${\rm Adm}(B_1)\subset {\rm Adm}(B)$, since
$K_{B_1}\subset K_B$. 
\smallskip

The following theorem answers this question.
To state it let us introduce the notion of a one-sided majorant,
which is natural, since the zeros are asymmetric.
We say that $w\in {\rm Adm}_+(\Theta)$ ($w\in {\rm Adm}_-(\Theta)$) if
there exists a nonzero function $f\in K_\Theta$ such that
$|f(x)|\le w(x)$ for a.e. $x>0$ 
(respectively for a.e. $x<0$ ).
\medskip
\\
{\bf Theorem 1.1.} {\it Let $w$ be a nonnegative function on
$[0,\infty)$ and let $\Omega=-\log w$. 

1. If  $w\in {\rm Adm}_+(B_1)$, then
\begin{equation}
\label{3}
\int\limits_1^\infty t^{-3/2} \Omega(t)dt<\infty.
\end{equation} 

2. If $w$ is 
positive, nonincreasing and the integral (\ref{3}) 
converges, then 
$w\in {\rm Adm}_+(B_1)$ and $w(|x|)\in {\rm Adm}(B_1)$.  }
\medskip

It turns out that there is more freedom in the behavior of the elements
of $K_{B_1}$ along the negative semiaxis (that is, when we are 
far from the poles). In particular, 
the majorants $w(t)=\exp(-A|t|^{1/2})$
are in ${\rm Adm}_-(B_1)$ (whereas, by Theorem 1.1,
the majorant of the form $\exp(-|t|^\alpha)$
belongs to ${\rm Adm}_+(B_1)$ if and only if $\alpha<1/2$).
Moreover, this result is sharp.
\medskip
\\
{\bf Theorem 1.2.} {\it The majorant $w(t)=\exp(-A|t|^{1/2})$ 
belongs to the class ${\rm Adm}_-(B_1)$ for any $A>0$. 
At the same time, if $|t|^{1/2}=o(\Omega(t))$
when $t\to -\infty$, then 
$w\notin {\rm Adm}_-(B_1)$. }
\medskip 

These results remain true for a wider class of Blaschke products
with zeros in a half-strip having certain density properties.
Let all $z_n$ lie in a half-strip $[0,\infty)\times [\delta, M]$,
where $M>\delta>0$. Assume that
\begin{equation}
\label{dens}
0<c\le \frac{{\rm card}\{n: {\rm Re}\, z_n\in [x,x+r]\}}{r}\le C<\infty,
\qquad r>r_0,\ \ x>0.
\end{equation}
\medskip
\\
{\bf Theorem 1.3.} {\it Let $B$ be a Blaschke product with zeros 
$z_n$ satisfying (\ref{dens}).
Then all statements of Theorems 1.1 and 1.2 are true for $B$ instead
of $B_1$.}
\medskip

Now we consider a class of Blaschke products with power growth
of zeros. Let $\beta>1/2$ and let 
$B_\beta$ be the Blaschke product with simple zeros at the points
$$
z_n=n^\beta+i, \quad n\in\mathbb{N}
$$
(if $\beta\le 1/2$, then the Blaschke condition does not hold).
Thus, the notation $B_1$ for the product in Theorems 1.1 and 1.2 agrees
with this new notation. We will study the asymptotic decay
of admissible majorants for $K_{B_\beta}$. Namely, let
$$
w_\alpha(x)=\exp(-|x|^\alpha), \qquad \alpha\in(0,1), 
$$
and put
$$
\alpha(\beta)=\sup\{\alpha: w_\alpha\in {\rm Adm}(B_\beta)\}.
$$
Analogously, one can define the numbers $\alpha_+(\beta)$ and
$\alpha_-(\beta)$. It follows from Theorems 1.1 and 1.2 that
$\alpha(1)=\alpha_+(1)=\alpha_-(1)=1/2$.
The following theorem shows that even such a rough characteristic
as $\alpha(\beta)$ has a rather complicated behavior.
\medskip
\\
{\bf Theorem 1.4.} 
\begin{equation}
\label{4}
\alpha(\beta)=\alpha_+(\beta)=
\begin{cases}
1/\beta, & \beta>2, \\
1/2,      & 2/3\le \beta\le 2, \\
-1+1/\beta, & 1/2<\beta<2/3;
\end{cases}
\end{equation}
\begin{equation}
\label{5}
\alpha_-(\beta)=
\begin{cases}
1/\beta, & \beta>2, \\
1/2,      & 1\le \beta\le 2, \\
1, & 1/2<\beta<1.
\end{cases}
\end{equation}
\smallskip
\\
{\bf Remarks.} 
1. It is natural to ask whether $w_{\alpha(\beta)}$ is admissible for
$K_{B_\beta}$. Our argument in Section 5 shows that for $\beta >2$ and 
some $\kappa>0$ the majorant $w_{\alpha(\beta)}^\kappa$ is admissible 
for $K_{B_\beta}$. If $2/3 <\beta \le 2$, then, by Theorem 3.1, 
$w_{\alpha(\beta)}^\kappa \notin {\rm Adm}(B_\beta)$ for any $\kappa>0$.
Our question remains open for $\beta\in (1/2,2/3)$.
\smallskip

2. An interesting feature of the limit exponent $\alpha(\beta)$
is that it is constant for $\beta\in[2/3,2]$, though for $\beta>1$
the zeros are sparse, whereas in the case $\beta<1$ the zeros
are much more dense (in particular, the sequence is not an interpolating one).
An analogous phenomenon may be observed in the problems of weighted polynomial 
approximation on discrete subsets of the line (see \S 5). 
These problems turn out to be 
closely related to admissibility conditions.
\smallskip

3. It should be noted that in the case $\beta<1$ the admissible
majorants on the negative semiaxis may decrease much faster
than on the positive one. However, 
any nonzero function $f\in K_{B_\beta}$, $2/3<\beta<1$, 
which is majorized on $(-\infty,0]$ by $w_\alpha$ with 
$\alpha\in (1/2,1)$, is unbounded on $[0,\infty)$ (see Remark 5.5).
\smallskip

4. Analogous results are obtained 
for the Blaschke products with two-sided
zeros having different power growth in positive and negative directions
(see Theorem 5.6).
We mention also that all these results may be easily generalized 
to the case of perturbed zeros with certain density properties.
\medskip

Finally we consider the case of zeros which approach the real
axis tangentially. Let $B$ be the Blaschke product with the  zeros
$z_n=n+iy_n$, where $0<y_n\le 1$, $n\in\mathbb{Z}$ 
(or $n\in\mathbb{N}$). The most interesting situation is when $y_n\to 0$, 
$|n|\to\infty$. In this case another phenomenon  occurs. If $y_n$ tend to zero 
not too rapidly, then there is no qualitative difference between 
the classes ${\rm Adm}(B)$ and ${\rm Adm}(e^{iaz})$, $a>0$. Otherwise, there are no
admissible majorants with more than power decay.
\medskip
\\
{\bf Theorem 1.5.} {\it Let the sequence $\{y_n\}_{n\in\mathbb{Z}}$
be even and nonincreasing for $n\ge 0$.
\par 
1. If 
\begin{equation}
\label{6}
\sum\limits_{n\in\mathbb{N}} n^{-2}\log\frac{1}{y_n}<\infty,
\end{equation}        
then any even majorant $w$, nonincreasing on $[0,\infty)$ and such that 
${\cal L}(w)<\infty$, is admissible for $K_B$. 
\smallskip
\par
2. Let $y:\mathbb{R}\to (0,\infty)$
be an even function nonincreasing on $[0,\infty)$
and such that $y(n)=y_n$, $n\in\mathbb{Z}$,
and let $Y=-\log y$.
If the function $Y(e^x)$ is convex on $\mathbb{R}$
and the series (\ref{6}) diverges, 
then any majorant $w$ 
such that, for any $N>0$, $w(t)=o(|t|^{-N})$, $t\to\infty$, 
is not admissible for $K_B$.}
\medskip

To obtain these results we combine the admissibility conditions 
and techniques of \cite{hm1,hm2} with some results on quasianalyticity.
\medskip

The paper is organized as follows. In \S 2 we present 
a few results of \cite{hm1,hm2} and some corollaries 
which we use later on.
In \S 3 we prove Theorem 1.1 whereas \S 4 is devoted to the proof 
of Theorems 1.2 and 1.3. 
Limit exponents for the Blaschke products with power growth 
of the zeros will be considered in \S 5. Finally, in \S 6 we discuss
the case of tangential zeros.

\bigskip

\begin{center}
{\bf \large \S 2. General sufficient conditions}
\end{center}

Here we present some of the results of the papers \cite{hm1, hm2, bh}.
We start with the following general criterion of admissibility. Here
$\Theta$ is an arbitrary, not necessarily meromorphic, inner function.
We denote by $\Pi$ the Poisson measure, that is, $d\Pi(t)=\frac{dt}{t^2+1}$.
If $\Omega=-\log w\in  L^1(\Pi)$, 
then there exists the Hilbert transform of the function
$\Omega$ defined as follows:
$$
\widetilde{\Omega}(x)=v.p.\, \frac{1}{\pi}
\int\limits_\mathbb{R} \left(\frac{1}{x-t} +\frac{t}{t^2+1}\right)\Omega(t)dt.
$$
\smallskip
\\
{\bf Theorem 2.1.} {\it A majorant $w$ with $\Omega \in  L^1(\Pi)$
is admissible  for $K_\Theta$ if and only if
there exists a function $m\in L^\infty(\mathbb{R})$ with $m\ge 0$, 
$mw\in L^2(\mathbb{R})$, and $\log m \in L^1(\Pi)$,
and an inner function $I$ such that 
\begin{equation}
\label{7}
\arg \Theta+2\widetilde{\Omega}=2\widetilde{\log m}+\arg I+2\pi k \qquad \mbox{a.e. on }\mathbb{R},
\end{equation}
where $k$ is a measurable function with integer values. Here $\arg\Theta$
is an arbitrary measurable branch of the argument.}
\medskip

A certain refinement of this parametrization formula for admissible 
majorants is obtained in \cite{bh}.
Namely, it is shown that the theorem remains true if we replace 
$\arg I$ in (\ref{7}) by a constant $\gamma\in\mathbb{R}$.
\medskip

The following theorem provides a condition sufficient for the 
representation of the form (\ref{7}). We denote by 
${\rm Osc} (f,I)$ the oscillation of a function 
$f$ on the set $I$, that is,
$$
{\rm Osc} (f,I)=\sup\limits_{s,t\in I}(f(s)-f(t)).
$$
\smallskip
\\
{\bf Theorem 2.2.} {\it Let $f$ be a $C^1$-function on $\mathbb{R}$ 
and let $\{d_n\}$ $($where $n\in\mathbb{Z}$ or $n\in\mathbb{N}$; 
in the latter case we assume that $d_1=-\infty$, and we do not require
$f$ to have a limit at infinity$)$ be an increasing 
sequence of real numbers such that $\lim_{|n|\to\infty} |d_n|=\infty$ and 
$$
f(d_{n+1})-f(d_n)\asymp 1,\qquad n\in\mathbb{Z}\ \  (n\ge2).
$$
Assume also that there is a constant $C>0$ such that
$$
{\rm Osc} (f,(d_n,d_{n+1}))\le C\qquad \mbox{and} \qquad
{\rm Osc} (f',(d_n,d_{n+1}))\le C
$$
for all $n\in\mathbb{Z}$ $(n\in\mathbb{N})$. Then $f$ admits the representation
$$
f=2\widetilde{\log m}+2\pi k+\gamma,
$$
where $m\in L^\infty(\mathbb{R})\cap L^2(\mathbb{R})$, $m\ge 0$,
$\log m \in L^1(\Pi)$, $\gamma \in\mathbb{R}$, and $k$ 
is a measurable integer-valued function. }
\medskip

This theorem is proved in \cite{hm2} under a small additional restriction
on the distances $d_{n+1}-d_n$ and in \cite{bh} in the general case. 
Such 
functions $f$ will be referred to as {\it mainly increasing functions}.
It should be mentioned that the condition ${\rm Osc} (f',(d_n,d_{n+1}))\le C$
may be replaced by a weaker integral estimate.
\medskip

Combining Theorems 2.1 and 2.2 we arrive at the following sufficient 
condition.
\medskip
\\
{\bf Corollary 2.3.} {\it Let $\Theta$ be a meromorphic inner function 
and let $\varphi$ be an increasing branch of the argument of $\Theta$.
If $\varphi+2\widetilde{\Omega}$ is a mainly increasing function, then $w\in{\rm Adm}(\Theta)$.}
\medskip

We have the following useful corollary.
\medskip
\\
{\bf Corollary 2.4.} {\it Let $\Theta_1$ and $\Theta_2$ 
be meromorphic inner functions with arguments $\varphi_1$ and 
$\varphi_2$ respectively. If the function $\varphi_1-\varphi_2$
is mainly increasing, then 
\smallskip
${\rm Adm}(\Theta_2)\subset{\rm Adm}(\Theta_1)$.}
\medskip
\\
{\bf Remark.} Clearly, under conditions of Corollary 2.4
analogous inclusions take place for the classes ${\rm Adm}_+(\Theta_j)$
and ${\rm Adm}_-(\Theta_j)$, $j=1,2$; namely, 
${\rm Adm}_+(\Theta_2)\subset {\rm Adm}_+(\Theta_1)$
and ${\rm Adm}_-(\Theta_2)\subset {\rm Adm}_-(\Theta_1)$.
Indeed, if $w\in {\rm Adm}_+(\Theta_2)$, then there exists a nonzero function
$f\in K_{\Theta_2}$ such that $|f(t)|\le w(t)$, $t<0$. Hence,
$|f|\in {\rm Adm}(\Theta_2)$ and, by Corollary 2.4, $|f|\in {\rm Adm}(\Theta_1)$.
Thus, $|f|\in {\rm Adm}_+(\Theta_1)$ and, consequently,
$w\in {\rm Adm}_+(\Theta_1)$.
\medskip

We get immediately the following corollary concerning
the inner functions with "almost linear"\, growth of the argument.
\medskip
\\
{\bf Corollary 2.5.} {\it Let $\Theta$ be a meromorphic inner function
such that $\varphi'\asymp 1$. Then there exist positive numbers $a$, $b$
such that ${\rm Adm}(e^{iaz})\subset{\rm Adm}(\Theta)\subset{\rm Adm}(e^{ibz})$.}
\smallskip
\\
{\bf Proof.}  Note that $\psi(t)=bt$ is 
the continuous argument of the inner function $e^{ibz}$.
Let $c\le \varphi'(t)\le C$, $t\in\mathbb{R}$, for some positive 
constants $c$, $C$. Take $b>C$. Then 
$bt-\varphi(t)$ is an increasing Lipschitz function
and, consequently, is mainly increasing.
Hence, by Corollary 2.4, 
${\rm Adm}(\Theta)\subset{\rm Adm}(e^{ibz})$. The proof of the 
second inclusion is analogous. $\bigcirc$
\medskip

Now we discuss another approach to admissible majorants
applicable to the case of Blaschke products with sparse zeros.
Let $\Theta$ be a meromorphic inner function with zeros $z_n$
repeated according to their multiplicities. Then there exists
an entire function $E$ in the Hermite--Biehler 
class (that is, $|E(z)| > |E(\overline {z})|$, $z\in {\mathbb{C}^+}$)
with zeros at the points $\overline z_n$
such that $\Theta=E^*/E$
(see, e.g., \cite[Lemma 2.1]{hm1}).
\smallskip
Here $E^* (z) = \overline {E (\overline z)}$.
                 
With the function $E$ we associate the de Branges space
${\cal H} (E) $ which consists of all entire functions 
$F$ such that $F/E$ and $F^*/E$ belong
to the Hardy class $H^2$.
The book of L. de Branges \cite{br} is devoted to the 
theory of spaces ${\cal H} (E)$; this theory  
has important applications in mathematical physics.
\smallskip

It is easy to see that the mapping $ F\mapsto F/E $ 
is a unitary operator from ${\cal H}(E)$ onto  $K_{\Theta_E}$,
where $\Theta_E=E^*/E$,
that is, $K_{\Theta_E}={\cal H}(E)/E$ (see, for example, \cite{bar99} or
\cite[Theorem 2.10]{hm1}).
\smallskip

An admissible majorant $w$ is said to be minimal if for any other 
admissible majorant $\tilde w$ such that $\tilde w\le Cw$ we have
$\tilde w \asymp w$, that is, $cw\le\tilde w\le Cw$ 
for some positive constant $c$. It turns out that the inclusion
$1\in{\cal H}(E)$ is crucial for the existence of a 
positive and continuous minimal majorant.
Namely, the following dichotomy is true (see \cite{bar04, bh, hm1}).
\medskip
\\
{\bf Theorem 2.6.} 
{\it Let $E$ be an entire function of zero exponential type such that
$|E(z)| > |E(\overline {z})|$, $z\in {\mathbb{C}^+}$.
Then, either
\par
a) $1/E\in L^2(\mathbb{R})$ and $1/|E|$ is 
(the unique up to equivalence) positive and continuous 
minimal majorant for $K_{\Theta_E} $;
\par 
b) $1/E\notin L^2(\mathbb{R})$ and there is no 
positive and continuous minimal majorant 
for $ K_{\Theta_E} $.}
\medskip

To conclude this section we show that 
${\rm Adm}(B)={\rm Adm}(e^{2\pi iz})$ if $B$ is the Blaschke
product with the zeros $z_n=n+i$, $n\in\mathbb{Z}$. Put $E_1(z)=e^{i\pi z}$
and $E_2(z)=\sin \pi(z+i)$. Then $B=E_2^*/E_2$. Clearly,
$|E_1(z)|\asymp |E_2(z)|$ for ${\rm Im}\, z\ge 0$ and so the spaces 
${\cal H}(E_1)$ and ${\cal H}(E_2)$ coincide as sets with 
equivalence of norms. Since 
$K_{E_1^*/E_1}={\cal H}(E_1)/E_1$ and $K_{E_2^*/E_2}={\cal H}(E_2)/E_2$
it follows that ${\rm Adm}(E_1^*/E_1)={\rm Adm}(E_2^*/E_2)$.

\bigskip

\begin{center}
{\bf \large \S 3. Proof of Theorem 1.1}
\end{center}

We start with the proof of Statement 1 of Theorem 1.1.
Moreover, we prove a somewhat stronger result. Recall that
$B_\beta$ is the Blaschke product with the zeros
$z_n=n^\beta+i$, $n\in\mathbb{N}$.
\medskip                   
\\
{\bf Theorem 3.1.} 
{\it Let $\beta>2/3$. If $f\in{K_{B_\beta}}$, then
$$
\int\limits_0^\infty \frac{\log |f(t)|}{t^{3/2}+1}dt>-\infty.
$$}

To prove Theorem 3.1
we make use of certain properties  of model subspaces. Recall that
the function
$$
k_z(\zeta)=
\frac{i}{2\pi}\cdot
\frac{1-\overline{\Theta(z)}\Theta(\zeta)}{\zeta-\overline z}
$$
is the reproducing kernel of the space ${K_{\Theta}}$
corresponding to a point $z\in\mathbb{C^+}$, that
is, $f(z)=(f,k_z)$, $f\in K_\Theta$, where $(\cdot,\cdot)$
stands for the usual inner product in $L^2(\mathbb{R})$.
In the case of meromorphic inner functions
the same formula gives the reproducing kernel at the point $z=x\in\mathbb{R}$.
\smallskip

Clearly, 
$$
|f(z)|\le\|f\|_2\|k_z\|_2,\qquad z\in\mathbb{C}^+ \cup \mathbb{R}
$$
(here we denote by $\mathbb{C}^+ \cup \mathbb{R}$ the closed upper half-plane).
Note also that 
$$
\|k_z\|_2^2=\frac{1-|\Theta(z)|^2}{4\pi{\rm Im}\, z},\quad z\in\mathbb{C}^+,\quad
\mbox{and}\quad
\|k_x\|_2^2=\frac{|\Theta'(x)|}{2\pi},\quad x\in\mathbb{R}.
$$
Now let 
$B=\prod_n b_n$ be a Blaschke product, where
$b_n(z)=e^{i\alpha_n}(z-z_n)/(z-\overline z_n)$. Clearly,
$1-|B(z)|^2 \le \sum_n (1-|b_n(z)|^2)$, $z\in\mathbb{C}^+$, and hence
$$
\frac{1-|B(z)|^2}{2{\rm Im}\, z}\le \sum\limits_n
\frac{1-|b_n(z)|^2}{2{\rm Im}\, z}=
\sum\limits_n \frac{2{\rm Im}\, z_n}
{|z-\overline z_n|^2}\le |B'({\rm Re}\, z)|, \quad z\in\mathbb{C}^+
$$
(see (\ref{1a})). Thus, we get the following estimate: if $f\in K_B$, then
\begin{equation}
\label{estim1}
|f(z)|\le C(f) |B'(x)|^{1/2},\qquad z=x+iy\in \mathbb{C}^+ \cup \mathbb{R}.
\end{equation}
 
We will also need to estimate the function $f\in K_B$
in the lower half-plane. It follows from the definition
of $K_B$ that the inclusion $f\in K_B$ implies that 
$f/B\in H^2(\mathbb{C}^-)$. Moreover, 
the function $g(z)=\overline{f(\overline z)}B(z)$, $z\in\mathbb{C}^+$, is in $K_B$.
Hence, applying the estimate (\ref{estim1}) we have
\begin{equation}
\label{estim2}
|f(z)|\le C(f)|B(z)|\cdot |B'(x)|^{1/2},\qquad z=x-iy\in\mathbb{C}^-
\end{equation}
(here the right-hand side may be infinite).
\smallskip

Consider the domain 
\begin{equation}
\label{dom}
\Delta=\mathbb{C}\setminus\{z:\, {\rm Re}\, z\ge 0, -2\le{\rm Im}\, z\le 0\}.
\end{equation}
Then each function $f\in{K_{B_\beta}}$ is analytic 
in $\Delta$. 
Let $\eta$ be a conformal mapping of the upper half-plane 
$\mathbb{C}^+$ onto the domain $\Delta$ such that $\eta(0)=0$, $\eta(\infty)=\infty$. 
By the Christoffel--Schwarz formula, $\eta$ is of the form
\begin{equation}
\label{conf}
\eta(z)=a_1+a_2\int\limits_{z_0}^{z}\zeta^{1/2}(\zeta-a)^{1/2}d\zeta
\end{equation}
where $a_1\in\mathbb{C}$, $a_2>0$, 
$z_0\in\mathbb{C}^+$ and $a>0$. Clearly,
$$
\eta(z)\sim a_2z^2/2,  \quad\mbox{ and }\quad 
\eta'(z)\sim a_2z
$$
when $|z|\to\infty$, $z\in\mathbb{C}^+ \cup \mathbb{R}$ (we write
$f(z)\sim g(z)$, $z\to z_0$, if 
\medskip
$\lim_{z\to z_0}f(z)/g(z)=1$).

In the proof of Theorem 3.1 we will use the following lemma.
\medskip                   
\\
{\bf Lemma 3.2.} 
{\it Let $f\in{K_{B_\beta}}$. If $\beta\ge 1$, then $f$ is bounded in $\Delta$.
If $1/2<\beta<1$, then
\begin{equation}
\label{estim3}
|f(z)|\le C_1\exp(C_2|z|^{-1+1/\beta}), \qquad z\in \Delta.
\end{equation}} 
{\bf Proof.} Since 
\begin{equation}
\label{bbeta}
|B'_\beta(x)|=\sum\limits_{n=1}^\infty
\frac{2}{(x-n^\beta)^2+1},
\end{equation}
it is easy to see that $B'_\beta\in L^\infty(\mathbb{R})$ for $\beta\ge 1$
and $B'_\beta\in L^\infty((-\infty,0])$ for $1/2<\beta<1$.
Let us show that for $1/2<\beta<1$
$$
|B'_\beta(x)|\asymp x^{-1+1/\beta},\qquad x>1.
$$
Let $x=t^\beta$. First we note that
$$
\sum\limits_{n=1}^{[t/2]}\frac{1}{(t^\beta-n^\beta)^2+1}+
\sum\limits_{n=[3t/2]}^\infty\frac{1}{(t^\beta-n^\beta)^2+1}
\le C t^{1-2\beta}=Cx^{-2+1/\beta}.
$$
We denote by $[s]$  the entire part of $s$. If
$[t/2]\le n\le [3t/2]$, then $|t^\beta-n^\beta|\asymp
|t-n|t^{\beta-1}$. Then
$$
\sum\limits_{[t/2]}^{[3t/2]}\frac{1}{(t^\beta-n^\beta)^2+1}
\asymp 
\sum\limits_{[t/2]}^{[3t/2]}\frac{1}{(t-n)^2t^{2\beta-2}+1}.
$$
Clearly,
$$
\sum\limits_{[t]-[t^{1-\beta}]}^{[t]}
\frac{1}{(t-n)^2t^{2\beta-2}+1}\le C t^{1-\beta}.
$$
On the other hand,
$$
\sum\limits_{[t/2]}^{[t]-[t^{1-\beta}]}
\frac{1}{(t-n)^2t^{2\beta-2}+1}
\asymp t^{2-2\beta}
\sum\limits_{[t/2]}^{[t]-[t^{1-\beta}]}\frac{1}{(t-n)^2}
\asymp t^{1-\beta}.
$$
The estimate of the sum over $[t]<n\le [3t/2]$ is analogous. 
\medskip

Now it follows from (\ref{estim1}) that 
$|f(z)|\le C$, $z\in\mathbb{C}^+ \cup \mathbb{R}$ in the case
$\beta\ge 1$ whereas 
for $1/2<\beta<1$ we have 
$|f(z)|\le C_1+C_2|z|^{-1/2 + 1/(2\beta)}$. To estimate $|f(z)|$ for 
$z\in\Delta\cap\mathbb{C}^-$, it suffices to estimate $|B_\beta(z)|$
from above.
\smallskip

Let $z=x-iy\in \Delta\cap\mathbb{C}^-$. Then
$$
2\log |B_\beta(z)|=
\sum\limits_{n=1}^\infty
\log\left(1+\frac{4y}{|z-\overline z_n|^2}\right)\le
4y\sum\limits_{n=1}^\infty \frac{1}{(x-n^\beta)^2+(y-1)^2}.
$$
By the estimates analogous to those above it is easily shown that
\begin{equation}
\label{l32}
\log |B_\beta(z)| \le
\begin{cases}
C, & \beta\ge 1, \\
C_1+C_2|z|^{-1+1/\beta}, & 1/2<\beta<1.
\end{cases}
\end{equation}
To complete the proof of the lemma we apply estimate (\ref{estim2}).
$\bigcirc$  
\medskip

We will also make use of the following simple lemma (see \cite[IIIG2]{ko1}).
\medskip                   
\\
{\bf Lemma 3.3.} 
{\it Let $g$ be a nonzero function analytic in $\mathbb{C}^+$ and continuous in
$\mathbb{C}^+ \cup \mathbb{R}$. If \linebreak$\log|g(z)|\le 
C|z|^\gamma$, $z\in \mathbb{C}^+ \cup \mathbb{R}$, 
where $0<\gamma<1$, then
$$
\int\limits_\mathbb{R}\frac{\log|g(t)|}{t^2+1}dt>-\infty.
$$ }
\\
{\bf Proof of Theorem 3.1.}
 Let $\beta>2/3$ and $f\in{K_{B_\beta}}$. By Lemma 3.2, $f$ is bounded in $\Delta$
if $\beta\ge 1$ and
$|f(z)|\le C_1 \exp(C_2|z|^{-1+1/\beta})$, $z\in \Delta$, if $2/3<\beta<1$.
\smallskip

Put $g(z)=f(\eta(z))$, $z\in\mathbb{C}^+$, where the conformal mapping
$\eta$ of $\mathbb{C}^+$ onto $\Delta$ is defined by (\ref{conf}).
Then $g \in H^\infty$ if $\beta\ge1$ and
$$
|g(z)|\le C_3\exp(C_4|z|^{-2+2/\beta}), \qquad z\in\mathbb{C}^+,
$$
if $2/3<\beta<1$.
Since $\beta>2/3$ we have $-2+2/\beta<1$. Hence, by Lemma 3.3,
$$
\int\limits_\mathbb{R} \frac{\log|g(t)|}{t^2+1}dt =
\int\limits_\mathbb{R} \frac{\log|f(\eta(t))|}{t^2+1}dt 
>-\infty.
$$
Taking into account that $ \eta(t)\asymp t^2$ and 
$\eta'(t)\asymp t$, $t\in\mathbb{R}$, $t\to+\infty$,
we obtain 
$$
\int\limits_0^\infty \frac{\log |f(t)|}{t^{3/2}+1} dt>-\infty.
\qquad \bigcirc
$$

Now we turn to the proof of Statement 2 (sufficiency)
of Theorem 1.1. In what follows we will make use of the following lemma.
\medskip                   
\\
{\bf Lemma 3.4.} 
{\it 
Let $\Omega=\log\frac{1}{w}$ be a nondecreasing 
function on $[0,\infty)$. If $f\in H^2$, $|f|\le 1$ 
a.e. on $\mathbb{R}$ and $|f|\le w$ a.e. on $(0,\infty)$, then 
$$
\log|f(z)|\le -\frac{1}{4}\Omega(|z|),\qquad z\in \mathbb{C}^+,\ \ {\rm Re}\, z>0.
$$ }
{\bf Proof.} Let $z=x+iy \in\mathbb{C}^+$ and $x>0$. By the Jensen inequality 
$$
\log|f(z)|\le\frac{y}{\pi}\int\limits_\mathbb{R} \frac{\log|f(t)|}{|t-z|^2}dt.
$$
Since $|f(t)|\le 1$ we have
$$
\log|f(z)|\le-\frac{y}{\pi}\int\limits_{x+y}^\infty \frac{\Omega(t)}{|t-z|^2}dt
\le - \frac{\Omega(|z|)}{\pi}\int\limits_{x+y}^\infty 
\frac{ydt}{(t-x)^2+y^2}dt=
-\frac{1}{4} \Omega(|z|). \qquad\bigcirc
$$

Recall that a sequence ${\lambda_n}\subset\mathbb{C}^+$ is said to be
an interpolating sequence if it satisfies the Carleson condition
$$
\inf\limits_n\prod\limits_{k\ne n}
\left|\frac{\lambda_n-\lambda_k}{\lambda_n-\overline \lambda_k}\right|
=\inf\limits_n 2{\rm Im}\,\lambda_n |B'(\lambda_n)|>0,
$$
where $B$ is the Blaschke product with the zeros $\lambda_n$.
By the Shapiro--Shields theorem (see \cite{nik,nk12}), 
in this case the rational fractions $k_n(z)=\frac{\sqrt{{\rm Im}\, \lambda_n}}
{z-\overline \lambda_n}$ form a Riesz basis in $K_B$ and, thus,
the inclusion $f\in K_B$ is equivalent to the representation 
\begin{equation}
\label{8}
f(z)=\sum\limits_n\frac{c_n\sqrt{{\rm Im}\, \lambda_n}}{z-\overline \lambda_n},
\qquad \{c_n\}\in\ell^2,
\end{equation}
with $\|f\|_2\asymp \|\{c_n\}\|_{\ell^2}$.
\smallskip

Clearly, the Blaschke product $B_1$ is interpolating. Hence,
each function $f$ in ${K_{B_1}}$ admits  the representation
\begin{equation}
\label{9}
f(z)=\sum\limits_{n\in\mathbb{N}} \frac{c_n}{z-n+i}.
\end{equation}

The following two lemmas 
relate the rate of decay of a function $f\in{K_{B_1}}$ with the 
properties of the coefficients $c_n$ in (\ref{9}).
Lemma 3.6 plays the key role in the proof of Statement 2 of Theorem 1.1.
However, for the sake of completeness we start with the
converse result.
\medskip                                     
\\
{\bf Lemma 3.5.} 
{\it Let $w$ be a positive nonincreasing
function on $[0,\infty)$, which tends to zero
faster than any power, that is,
\begin{equation}
\label{10}
\lim\limits_{t\to\infty}t^N w(t)=0
\end{equation}
for any $N>0$. If $w\in {\rm Adm}_+(B_1)$, then there exists a nonzero sequence 
$\{c_n\}_{n\in\mathbb{N}}$ such that
\begin{equation}
\label{11}
\log |c_n|\le -C\Omega(n)
\end{equation}
and all the moments of the sequence $\{c_n\}$ are equal 
to zero, that is,
\begin{equation}
\label{12}
\sum\limits_{n=1}^\infty c_n n^k=0,\qquad k\in \mathbb{Z}_+.
\end{equation} }

In the converse statement we impose certain regularity conditions
on the majorant $w$. We say that a majorant $w$ is {\it regular} if
$w$ is even, $0<w\le 1$, $w$ is nonincreasing on $[0,\infty)$ 
and the function $t\Omega'(t)$ is nondecreasing on $[0,\infty)$
(recall that $\Omega=\log\frac{1}{w}$).
The last property is equivalent to the following:
the function $G(s)=\Omega(e^s)$ is \medskip
convex.
\\
{\bf Lemma 3.6.} 
{\it Let $w$ be an even regular majorant on $\mathbb{R}$,
which tends to zero faster than any power.
If there exists a nonzero sequence $\{c_n\}_{n\in\mathbb{N}}$ 
such that
$$
|c_n|\le n^{-3}w(n), \qquad n\in \mathbb{N},
$$
and equalities (\ref{12}) hold, then $w\in{\rm Adm}(B_1)$. }
\medskip
\\
{\bf Proof of Lemma 3.5.}
Let $f$ be a function of the form (\ref{9}) such that
$|f(t)|\le w(t)$, $t>0$. By Lemma 3.2,
${K_{B_1}}\subset L^\infty(\mathbb{R})$. Thus, without loss of generality we may 
assume that $|f(t)|\le 1$, $t\in\mathbb{R}$.
\smallskip
                     
Since $f\in{K_{B_1}}$, it follows that the function $g(z)=f(z)/B(z)$, $z\in\mathbb{C}^-$,
is in $H^2(\mathbb{C}^-)$. We have $|g(t)|=|f(t)|$, $t\in\mathbb{R}$,
and $|c_n|=|g(\overline z_n)|/|B'(z_n)|$. Applying Lemma 3.4
to $g$ we get the estimate 
$\log |g(\overline z_n)|\le -C\Omega(n)$, which implies (\ref{11}).
\smallskip

Now we will show that (\ref{12}) is fulfilled.
By hypotheses, $|f(t)|=o(t^{-N})$, $t\to+\infty$, for any $N>0$.
By Lemma 3.2, $f$ is bounded in the domain $\Delta$ defined by
(\ref{dom}), and therefore we also have 
$|f(t)|=o(|t|^{-N})$, $t\to-\infty$, for any $N>0$.
Since $f$ tends to zero faster than any power, the function $(z+i)^kf(z)$ 
is in the Hardy class $H^1(\mathbb{C}^+)$ for every $k\ge 0$. Hence,
$$
\int\limits_{\mathbb R}(t+i)^k\left(\sum\limits_{n=1}^\infty
\frac{c_n}{t-n+i}\right)dt=0.
$$
Therefore, for $N>k+1$,
$$
\lim_{A\to\infty}\int\limits_{\mathbb R}\left( \frac{Ai}{Ai-t} \right)^N
(t+i)^k\left(\sum\limits_{n=1}^\infty\frac{c_n}{t-n+i}\right)dt
$$
$$
= \lim_{A\to\infty} 
\sum\limits_{n=1}^\infty c_n
\int\limits_{\mathbb R}
\left( \frac{Ai}{Ai-t} \right)^N 
\frac{(t+i)^k}{t-n+i}\, dt =0.
$$
By the residue calculus (in the lower half-plane), we have
$$
\lim_{A\to\infty}\sum\limits_{n=1}^\infty  
\left(\frac{Ai}{Ai-n+i}\right)^N n^kc_n=0
$$
and, finally,
$$
\sum\limits_{n=1}^\infty n^kc_n=0. \qquad \bigcirc
$$
\smallskip
\\
{\bf Proof of Lemma 3.6.}
Let us define $f \in K_{B_1}$ by formula (\ref{9}). Then (\ref{12})
implies for any $k\in\mathbb{N}$
$$
(t+i)^kf(t)-\sum\limits_{n\in\mathbb{N}} \frac{c_n n^k}{t-n+i}=0.
$$
Hence,
$$
|f(t)|\le\inf\limits_{k\in\mathbb{Z}_+}\frac{1}{|t|^k}\sum\limits_{n\in\mathbb{N}}
|c_n|n^k.
$$                 
Since $|c_n|\le n^{-3}w(n)$, we have
$$
|f(t)|\le C\inf\limits_{k\in\mathbb{Z}_+}\frac{1}{|t|^k}\sup\limits_{n\in\mathbb{N}}
[w(n)n^{k-1}].
$$
It is easy to show that for $|t|\ge 1$
$$
\inf\limits_{k\in\mathbb{Z}_+}\frac{1}{|t|^k}\sup\limits_{n\in\mathbb{N}}
[w(n)n^{k-1}] \le
\inf\limits_{p\ge 0}\frac{1}{|t|^p}\sup\limits_{q\ge 1}
[w(q)q^p].
$$
Recall that if $G$ is a function on $[0,\infty)$, then
its Legendre transform $G^\#$ is defined as
$$
G^\#(x)=\sup\limits_{t\ge 0}(xt-G(t)).
$$
Put $G(t)=\Omega(e^t)$. Then
$$
\sup\limits_{q\ge 1}w(q)q^p=\sup\limits_{q\ge 1}e^{p\log q-\Omega(q)}
=e^{G^\#(p)}
$$
and 
$$
\inf\limits_{p\ge 0}e^{G^\#(p)-p\log |t|}=e^{-(G^\#)^\#(\log |t|)}.
$$
It is well known that for convex functions $G$ we have
$(G^\#)^\#=G$ and so
$$
|f(t)|\le e^{-G(\log |t|)}=w(t).
$$
The proof is completed. $\bigcirc$
\medskip

Now we will use the following theorem on analytic
quasianalyticity, which is due to P. Koosis \cite{ko01}.
Many results of this type were obtained by 
M.M. Dzhrbashyan, L. Carleson, B. R.-Salinas and B.I. Korenblum 
in 1940-s--1960-s.
We denote by $C_A^\infty$ the class of functions
analytic in the unit disc $\mathbb{D}$ and infinitely differentiable in its
closure.
\medskip
\\
{\bf Theorem 3.7.} {\it Let $\{w(n)\}_{n\in\mathbb{N}}$ be a positive sequence 
such that $w(n)=o(n^{-l})$, $n\to\infty$, for any $l>0$.
Put for $z\in\mathbb{C}$
$$
w_*(z)=\sup\{p(z):\  p \mbox{  is a polynomial and  }\ 
|p(n)w(n)|\le 1, n\in\mathbb{Z}_+\}.
$$
Then the following statements are equivalent:

1. There exists a nonzero function $f\in C_A^\infty$, 
$f(z)=\sum_{n\ge 0}c_nz^n$, such that $|c_n|\le w(n)$
and $f^{(k)}(1)=0$, $k\in\mathbb{Z}_+$;

2. 
$$
\sum\limits_{n=1}^\infty n^{-3/2}\log w_*(n)<\infty.
$$}
 
We will use only implication 2$\Longrightarrow$1.
Note that $w_*(n)\le 1/w(n)$. Thus,  if the integral (\ref{3})
converges, then there exists a nonzero function $f$ as 
in Statement 1 above.
\medskip

Now Theorem 1.1 follows almost immediately from Theorems 3.1 and 3.7.
\medskip
\\
{\bf Proof of Theorem 1.1.} 
Statement 1 is a particular case of Theorem 3.1.
\smallskip

Assume that $w$ is a positive nonincreasing function 
such that the integral
(\ref{3}) converges. Without loss of generality we
may assume that $w(t)\equiv 1$, $t\in [0,1]$. 
We replace $w$ by a smaller
regular majorant $w_1$ (see the definition given before Lemma 3.6). 
Put
$$
\Omega_1(x)=\int\limits_{0}^{ex}\frac{\Omega(t)}{t}dt.
$$
Then
$$
\Omega_1(x)\ge\int\limits_{|x|}^{e|x|}\frac{\Omega(t)}{t}dt\ge \Omega(|x|)
$$
and, consequently, $w_1=e^{-\Omega_1}\le e^{-\Omega}=w$.
Since $\Omega_1'(x)=\Omega(ex)/x$, the function $x\Omega_1'(x)$
is nondecreasing and so $w_1$ is regular.
Finally we put $w_2(t)=(t+1)^{-3}w_1(t)$, $t\ge 0$.
It is clear that the convergence of the integral (\ref{3})
implies that
$$
\int\limits_1^\infty t^{-3/2} \log\frac{1}{w_2(t)} dt<\infty.
$$
Then, by Theorem 3.7,
there exists a nonzero function $f\in C_A^\infty$, $f(z)=
\sum_{n=1}^\infty c_n z^n$,
such that $|c_n|\le w_2(n)$, $n\in\mathbb{N}$, 
and $f^{(k)}(1)=0$, $k\in\mathbb{Z}_+$. The latter condition means that
$$
\sum\limits_{n=1}^\infty c_n n^k=0,\qquad k\in\mathbb{Z}_+.
$$

Since $w_1$ is a regular majorant and $|c_n|\le n^{-3}w_1(n)$, $n\in\mathbb{N}$,
it follows from Lemma 3.6 that $w_1\in{\rm Adm}(B_1)$.
Hence, $w(|x|)\in{\rm Adm}(B_1)$ and, in particular,
$w\in{\rm Adm}_+(B_1)$. $\bigcirc$
         
\bigskip

\begin{center}
{\bf \large \S 4. Majorants on the negative semiaxis}
\end{center}
{\bf Proof of Theorem 1.2.} We start with the proof of the sharpness
of the exponent $1/2$ on the negative semiaxis.
Assume that 
\begin{equation}
\label{14a}
|t|^{1/2}=o(\Omega(t)), \quad t\to -\infty.
\end{equation}
We will show that $w\notin {\rm Adm}_-(B_1)$.
\smallskip

Let $\Delta$ be the domain defined by (\ref{dom}).
Then, by Lemma 3.2,
each function $f\in{K_{B_1}}$ is analytic and bounded in $\Delta$.  
Recall that we denote by $\eta$ the conformal mapping (\ref{conf})
of $\mathbb{C}^+$ onto $\Delta$ and we have  
\begin{equation}
\label{14d}
\eta(z)\sim a_2z^2/2,  \qquad |z|\to\infty, \ \  z\in\mathbb{C}^+,
\end{equation}
where $a_2>0$. Hence, $\eta^{-1}(-x)\in \{z: {\rm Im}\, z\ge |{\rm Re}\, z|\}$
for sufficiently large $x>0$. Also, if we put 
$\Gamma= \eta^{-1}((-\infty,0])$, then for $z\in\Gamma$
we have $|\eta(z)|\ge c|z|^2$ for some $c>0$ when $|z|$ is sufficiently large.
\smallskip

Let $f\in{K_{B_1}}$ and $|f(t)|\le w(t)$, $t<0$. Put $g(z)=f(\eta(z))$.
Then $g$ is a bounded analytic function in $\mathbb{C}^+$.
By (\ref{14a}), for $z\in\Gamma$ we have
\begin{equation}
\label{14e}
|g(z)|\le w(|\eta(z)|)\le e^{-C|z|}, \qquad |z| \to\infty,
\end{equation}
for any $C>0$. The estimate (\ref{14d}) and an elementary argument using
the theorem on two constants permit us to obtain (\ref{14e}) for any $C>0$ and
$z\in i\mathbb{R}_+$. Hence, by the Phragm\'{e}n--Lindel\"of
principle, $g\equiv 0$ in $\mathbb{C}^+$.
\medskip

Now we turn to the proof of admissibility of 
the majorants $w(t)=\exp(-A|t|^{1/2})$ 
on the negative semiaxis. In fact, we prove the following
 stronger result: 
\medskip
\\
{\bf Theorem 4.1.} {\it Let $B_\beta$ be the Blaschke product 
with the zeros  $n^\beta+i$, $n\in\mathbb{N}$, $\beta>1/2$, and let
$$
W_A(t)=\exp(-A|t|^{1/2}).
$$

1. If $A<\pi$, then $W_A \in {\rm Adm}_-(B_2)$.
\smallskip

2. If $1\le \beta<2$, then the majorant $W_A$ belongs to the class 
${\rm Adm}_-(B_\beta)$ for any 
\medskip
$A>0$.  }
\\
{\bf Proof of Statement 1.}
We consider an auxiliary 
Blaschke product $B^\circ$
with the zeros at the points $-1+i$, $i$ and
$(\rho n)^2+i$, $n\in\mathbb{N}$. We will show that 
$W_A\in{\rm Adm}_-(B^\circ)$ for any $A<\pi/\rho$.
\smallskip

Consider the entire function
\begin{equation}
\label{15b}
E(z)=\prod\limits_{n\in\mathbb{N}}\left(1-\frac{z}{(\rho n)^2-i}\right)
= c\, \frac{\sin (\pi\rho^{-1}\sqrt{z+i})}{\sqrt{z+i}}
\end{equation}
where $c$ is some absolute constant. Then we have
\begin{equation}
\label{15}
\log |E(x)|\sim \frac{\pi|x|^{1/2}}{\rho}, \qquad x\to-\infty.
\end{equation}
Now let $ x\in (k-1/2,k+1/2)$, $k\in\mathbb{N}$. We write
$$
|E(\rho^2 x^2)|= \left| \frac{\rho^2 k^2-\rho^2x^2-i}{\rho^2k^2+i} \right|
\prod\limits_{n\ne k} \left| 1-\frac{\rho^2 x^2}{\rho^2 n^2-i} \right|.
$$
It is easy to see that there exist positive constants $C_1$ and $C_2$
depending on $\rho$ but not depending on $k$ and $x$ such that
$$
C_1 \le \prod\limits_{n\ne k} 
\left| 1-\frac{\rho^2 x^2}{\rho^2 n^2} \right|^{-1}
\left| 1-\frac{\rho^2 x^2}{\rho^2 n^2-i} \right| \le C_2.
$$
Hence
$$
|E(\rho^2 x^2)| \asymp
\left| \frac{\rho^2 k^2-\rho^2x^2-i}{\rho^2k^2+i} \right|
\cdot \frac{| \sin\pi x|}
{\pi |x| } \left|1-\frac{x^2}{k^2}\right|^{-1},
$$
and it follows that
\begin{equation}
\label{15a}
\frac{C_3} {k^2} \le |E(\rho^2 x^2)| \le \frac{C_4} {k}.
\end{equation}
Thus, 
\begin{equation}
\label{16}
|E(x)|\ge C_5 x^{-1}, \qquad x>1. 
\end{equation}

Now, put ${E^\circ}(z)=(z+i)(z+1+i)E(z)$. 
Then $B^\circ=(E^\circ)^*/E^\circ$.
It follows from
the estimates (\ref{15}) and (\ref{16}) that for any $A<\pi/\rho$
there exist positive constants $C_6$ and $C_7$ such that
$$
|E^\circ(x)|^{-1}\le C_6\exp(-A|x|^{1/2}), \qquad x<0,
$$
and
$$
|E^\circ(x)|^{-1}\le C_7|x|^{-1}, \qquad x>1.
$$

Hence, $1/E^\circ \in L^2(\mathbb{R})$. It follows immediately
from the definition of $E^\circ$ that $|E^\circ(x+iy)|$
is an increasing function of $y\ge 0$ for each $x \in\mathbb{R}$
and so $1/E^\circ \in H^2$. Clearly, in this case $1/E^\circ \in 
K_{B^\circ}$ and, in particular,
$1/|E^\circ| \in {\rm Adm}_-(B^\circ)$
(by Theorem 2.6,  
the function $1/|E^\circ|$ is a positive minimal admissible majorant
for $K_{B^\circ}$).
In particular, the majorant $W_A$ is in ${\rm Adm}_-(B^\circ)$ for any 
\smallskip
$A<\pi/\rho$. 

Now let $\rho>1$.
To prove Statement 1 
it suffices to show that ${\rm Adm}_-(B^\circ)\subset {\rm Adm}_-(B_2)$.
Let $\varphi$ be an increasing branch of the argument of the Blaschke product
$B_2$ and let $\psi$ be an increasing branch of the argument of 
$B^\circ$. Then
$$
\varphi'(t)=2\sum\limits_{n\in\mathbb{N}}\frac{1}{(t-n^2)^2+1},
$$
$$
\psi'(t)=\frac{2}{t^2+1}+\frac{2}{(t+1)^2+1}
+2\sum\limits_{n\in\mathbb{N}}\frac{1}{(t-\rho n^2)^2+1}.
$$

For $n\in\mathbb{N}$ and $M>0$ put $d_n=(Mn)^2$. 
Since $\rho>1$, for sufficiently large $M$ we have
$$
\int\limits_{d_n}^{d_{n+1}}\varphi'(t)dt\asymp 1
\qquad \mbox{and}\qquad
\int\limits_{d_n}^{d_{n+1}}(\varphi'(t)-\psi'(t))dt\asymp 1.
$$
Note also that there is $C>0$ such that the function $\varphi-\psi$
is an increasing Lipschitz function for $t<-C$.
Hence, $\varphi-\psi$ is mainly increasing and, by the remark
after Corollary 2.4, ${\rm Adm}_-(B^\circ)\subset {\rm Adm}_-(B_2)$. 
Thus, the majorant $W_A$
is in ${\rm Adm}_-(B_2)$ for any $A<\pi/\rho$. Since $\rho$
is an arbitrary number greater than $1$, the proof of Statement 1
is completed.
\medskip
\\
{\bf Proof of Statement 2.}
Let $B^\circ$ be the same Blaschke product
with the zeros $-1+i$, $i$ and
$(\rho n)^2+i$, $n\in\mathbb{N}$. But this time we will assume
$\rho$ to be small. We have shown in the proof
of Statement 1 that $W_A$ belongs to ${\rm Adm}_-(B^\circ)$
for any $A< \pi/\rho$. 
\smallskip

Recall that we denote by $\varphi_\beta$ 
an increasing continuous branch of the argument of $B_\beta$
(see formula (\ref{bbeta})).
If $1\le \beta <2$, then it is easy to see that 
$\varphi_\beta-\psi$ is a mainly increasing function 
for any $\rho>0$ (take $d_n= (Mn)^{\beta}$, $n\in\mathbb{N}$,
for a sufficiently large $M$). 
Now, by Corollary 2.4,
${\rm Adm}(B^\circ)\subset {\rm Adm}(B_\beta)$ for any $\rho>0$
and therefore $W_A \in{\rm Adm}_-(B_\beta)$ for any $A>0$.
$\bigcirc$
\medskip
\\
{\bf Remark.} It is easy to see that the constant $\pi$ in Statement 1
is sharp, that is, $W_A\notin {\rm Adm}_-(B_2)$ for
$A\ge \pi$. Indeed, let $E$ be the function defined in 
(\ref{15b}) with $\rho=1$.
Then $|E(x)|\asymp |x|^{-1/2} \exp(\pi|x|^{1/2})$, $x<-1$,
and, by (\ref{15a}), $|E(x)|\le C x^{-1/2}$, $x>1$.
Assume that $f\in K_{B_2}$ and $|f(t)|\le W_A(t)$, $t<0$, where $A>\pi$.
Note that, by Lemma 3.2, $f$ is bounded on $\mathbb{R}$.
It is easily seen that $F=fE$ is an entire function of order at most $1/2$ 
(see \cite[Theorem 3.1]{hm1}), which is
bounded on the real axis and $|F(t)|\to 0$, $t\to-\infty$. Hence, $F\equiv 0$.
\medskip

The following statement shows that the fast decay on the
negative semiaxis is compatible with any admissible decay 
on the positive semiaxis from Theorem 3.1.
\medskip
\\
{\bf Corollary 4.2.} {\it For any $A>0$ and for any nonincreasing majorant 
$w$ with the finite  integral (\ref{3}) 
there exists a nonzero function $f\in K_{B_1}$ such that
$$
|f(t)|\le w(t),\quad t>0, \quad\mbox{and}\quad |f(t)|\le W_A(t),\quad t<0.
$$ }
{\bf Proof.}
Let $\widetilde{B}$ be the Blaschke product with zeros $\tau n+i$, $n\in\mathbb{N}$,
where $\tau>1$. Clearly, Theorem 1.1 is applicable also to 
the space $K_{\widetilde{B}}$ and, in particular, 
each nonincreasing majorant $w$ such that the integral (\ref{3})
converges  is in the class 
\smallskip
${\rm Adm}_+(\widetilde{B})$. 

Now let $D=B_{3/2}\widetilde{B}$, where $B_{3/2}$
is the Blaschke product with the zeros $n^{3/2}+i$, $n\in\mathbb{N}$. 
Note that, by Lemma 3.2, $K_{B_{3/2}}\subset L^\infty(\mathbb{R})$
and $K_{\widetilde{B}}\subset L^\infty(\mathbb{R})$.
Therefore, if $f\in K_{B_{3/2}}$ and $g\in K_{\widetilde{B}}$, 
then $fg\in K_D$.
\smallskip

By Statement 2 of Theorem 4.1,
there exists $f\in K_{B_{3/2}}$ such that $|f(t)|\le W_A(t)$,
$t<0$, and there exists $g\in K_{\widetilde{B}}$ such that $|g(t)|\le w(t)$, 
$t>0$. Hence the majorant 
$$
W(t)=
\begin{cases}
W_A(t), & t<0, \\
w(t),   & t>0. 
\end{cases}
$$
is admissible for $K_D$, since $|f(t)g(t)|\le CW(t)$, $t\in\mathbb{R}$.
\smallskip

To complete the proof note that $\arg B_1- \arg D$ is a mainly
increasing function since $\tau>1$. Now we apply Corollary 2.4
and see that $W\in{\rm Adm}(B_1)$. $\bigcirc$
\medskip

We conclude this section with the proof of Theorem 1.3.
\medskip
\\
{\bf  Proof of Theorem 1.3.}
For $h,\rho>0$ let $B_{h,\rho}$ be the Blaschke product with the zeros
$n/\rho +ih$, $n\in\mathbb{N}$. It is clear that all results 
of Theorems 1.1 and 1.2 are true also for the products $B_{h,\rho}$.
Denote by $\varphi_{h,\rho}$ an increasing branch of the argument
of $B_{h,\rho}$. 
\smallskip

Let all $z_n$ lie in a half-strip $[0,\infty)\times [\delta, M]$ 
and satisfy the condition (\ref{dens}).
Then there exist $L>0$ and $N\in\mathbb{N}$ such that 
each rectangle $[x,x+L]\times [\delta, M]$ contains at least one and 
no more than $N$ of the zeros $z_n$. 
It is easy to see that $\varphi'(x)\asymp 1$, $x>0$ \cite[Theorem 3.4]{hm2},
and $\varphi'(x)\asymp {|x|^{-1}}$, $x<-1$,
where $\varphi$ is an increasing branch of the argument of $B$.
Then there exist positive numbers $h_1$, $\rho_1$, $h_2$, and $\rho_2$
such that
$$
\varphi'_{h_1,\rho_1}(x)< \varphi'(x)< \varphi'_{h_2,\rho_2}(x),\qquad x\in\mathbb{R},
$$
and, thus, the functions $\varphi-\varphi_{h_1,\rho_1}$
and $\varphi_{h_2,\rho_2}-\varphi$ are increasing. By Corollary 2.4,
$$
{\rm Adm}(B_{h_1,\rho_1})\subset{\rm Adm}(B) \subset {\rm Adm}(B_{h_2,\rho_2}),
$$
$$
{\rm Adm}_+(B_{h_1,\rho_1})\subset{\rm Adm}_+(B)\subset {\rm Adm}_+(B_{h_2,\rho_2}),
$$
and 
$$
{\rm Adm}_-(B_{h_1,\rho_1})\subset{\rm Adm}_-(B)\subset {\rm Adm}_-(B_{h_2,\rho_2}).
\qquad \bigcirc     
$$

\bigskip

\begin{center}
{\bf \large \S 5. Power growth of zeros}
\end{center}

The proof of Theorem 1.4 consists of a few steps.
First we obtain the formulas (\ref{4})  and (\ref{5}) 
for the case $\beta\ge 1$.
Then, making use of a method of \cite{hm2}
we complete the proof of the formula (\ref{4}). 
Finally, we will show that $\alpha_-(\beta)=1$ 
for 
\smallskip
$\beta<1$.

We will use repeatedly the following lemma.
\medskip
\\
{\bf Lemma 5.1.} {\it Let $\beta>\gamma>1/2$. Then
${\rm Adm}(B_\beta)\subset{\rm Adm}(B_\gamma)$. 
In particular, the functions $\alpha$, $\alpha_+$ and $\alpha_-$
are nonincreasing functions of $\beta>1/2$. }
\smallskip
\\
{\bf Proof.}
Denote by $\varphi_\beta$ an increasing branch of the argument of 
$B_\beta$.  Then
$$
\varphi'_\beta(t)=2\sum\limits_{n\in\mathbb{N}}\frac{1}{(t-n^\beta)^2+1}.
$$
It is easy to see that for $\beta\ge 1$ we have
\begin{equation}
\label{17}
\varphi'_\beta(t)\asymp \frac{1}{(t-(n_\beta(t))^\beta)^2+1}, \qquad t>0,
\end{equation}
where $n_\beta(t)$ is the integer closest to $t^{1/\beta}$.
On the other hand, for $1/2<\beta<1$ we have shown in the proof of 
Lemma 3.2 that
\begin{equation}
\label{17a}
\varphi'_\beta(t)\asymp t^{-1+1/\beta}, \qquad t>1.
\end{equation}
Finally, for $\beta>1/2$,
\begin{equation}
\label{17b}
\varphi'_\beta(t)\asymp |t|^{-2+1/\beta}, \qquad t<-1.
\end{equation} 
\smallskip

Let $\beta>\gamma$. For $M>0$ and $n\in\mathbb{N}$ put 
$d_n= (Mn)^\gamma$. It follows from 
(\ref{17}) and (\ref{17a}) that for sufficiently large $M$ we have
$$
\int\limits_{d_n}^{d_{n+1}}\varphi'_\gamma(t)dt\asymp 1
\qquad\mbox{and}\qquad
\int\limits_{d_n}^{d_{n+1}}(\varphi'_\gamma(t)-\varphi'_\beta(t))dt\asymp 1.
$$
It is also clear that the function 
$\varphi_{\gamma}-\varphi_\beta$ is increasing when $t<0$
and $|t|$ is sufficiently large.
Thus, the function $\varphi_{\gamma}-\varphi_\beta$ is mainly increasing.
By Corollary 2.4, ${\rm Adm}(B_\beta)\subset {\rm Adm}(B_\gamma)$ and, consequently,
$\alpha(\beta)\le\alpha(\gamma)$. Analogously, by the remark 
following Corollary 2.4,
$\alpha_+(\beta)$ and $\alpha_-(\beta)$ are nonincreasing functions of 
$\beta$. $\bigcirc$
\medskip

We will also need the following lemma on asymptotics of certain
canonical products.
\medskip
\\
{\bf Lemma 5.2.} {\it For $\beta>2$ consider the entire function 
\begin{equation}
\label{18}
E_\beta(z)=\prod\limits_{n=1}^\infty
\left(1-\frac{z}{n^\beta-i} \right).
\end{equation} 
Then 
$$
\log|E_\beta(x)|\asymp |x|^{1/\beta},\qquad |x|\to\infty.
$$}
\\
{\bf Proof.} The asymptotics of $E_\beta$
outside an exceptional set is given in 
\cite[Chapter II, Theorem 5]{lev}. 
For $x\to-\infty$ we have
$\log|E_\beta(x)|\asymp |x|^{1/\beta}$.
For sufficiently small $\varepsilon>0$, $\varepsilon_1>0$ we have the estimate
$$
\log |E_\beta(z)|\asymp |z|^{1/\beta}, \qquad
|z|\to\infty, \ \ |\arg z|<\varepsilon, \ \ z\notin  \cup_n D_n,
$$
where $D_n=\{w\in\mathbb{C}: |w-(n^\beta-i)|<\varepsilon_1 n^{\beta-1}\}$.
Dividing $E_\beta$ by $(z-(n^\beta-i))$ and applying the maximum principle
in the discs $D_n$ we conclude that 
$\log |E_\beta(x)|\asymp x^{1/\beta}$, $x\to\infty$. $\bigcirc$
\medskip
\\
{\bf Proof of (\ref{4}) and (\ref{5}) 
for the case $\beta\ge 1$.}
First, assume that $\beta>2$. 
Clearly, $B_\beta=E_\beta^*/E_\beta$, where the entire function
$E_\beta$ is defined by (\ref{18}). By Lemma 5.2,
$\log|E_\beta(x)|\asymp |x|^{1/\beta}$, $|x|\to\infty$.
Hence, $1\in{\cal H}(E_\beta)$ and, by Theorem 2.6, 
$1/|E_\beta|$ is the minimal admissible majorant for $K_{B_\beta}$.
Thus, $\alpha(\beta)=1/\beta$. 
\smallskip

Also we have $\alpha_+(\beta)\ge1/\beta$ and $\alpha_-(\beta)\ge 1/\beta$.
To prove the converse inequalities let us show that
$w_\gamma$ does not belong to ${\rm Adm}_+(B_\beta)$ or ${\rm Adm}_-(B_\beta)$ 
when $\gamma>1/\beta$. Indeed, if a function $f\in K_{B_\beta}$ 
satisfies $|f(t)|\le e^{-t^\gamma}$, $t>0$,
then $F=fE_\beta$ is an entire function of order at most $1/\beta<1/2$ and
$|F(t)|\to 0$, $t\to +\infty$. Therefore, 
by the  Phragm\'en--Lindel\"of principle, $F\equiv 0$.
The same argument works for $\alpha_-(\beta)$.
\smallskip

So, we have shown that for $\beta>2$ 
$$
\alpha(\beta)=\alpha_+(\beta)=\alpha_-(\beta)=1/\beta.
$$ 

By Lemma 5.1, the functions $\alpha$, $\alpha_+$ and $\alpha_-$
are nonincreasing. Since, by Theorem 1.1, $\alpha(1)=1/2$ and 
$\lim_{\beta\to 2+0}\alpha(\beta)=1/2$,
we see that $\alpha(\beta)=1/2$ for $1<\beta\le 2$.
In a similar way we get 
$\alpha_+(\beta)=\alpha_-(\beta)=1/2$ for $1<\beta\le 2$. $\bigcirc$
\medskip

Now we obtain an estimate from below for $\alpha(\beta)$
in the case $\beta<1$.
\medskip
\\
{\bf Lemma 5.3.}
{\it Let $1/2<\beta<1$. Let $w$ be an even function nonincreasing on 
$[0,\infty)$ and satisfying the condition
$$
\Omega(t)= o(t^{-1+1/\beta}),\qquad t\to\infty.
$$
Then $w\in{\rm Adm}(B_\beta)$.}
\smallskip
\\
{\bf Proof.} We apply a method of the paper \cite{hm2}, which was used there
to deduce the admissibility of even majorants with convergent 
logarithmic integral from Corollary 2.3 in the case $\varphi'\asymp 1$.
\smallskip

Without loss of generality let $\Omega(t)=0$, $|t|\le 1$. We regularize
the majorant $w$ by considering the majorants
$$
\Omega_1(x)=\int\limits_{0}^{e|x|}\frac{\Omega(t)}{t}dt
$$
and 
$$
\Omega_2(x)=\int\limits_{0}^{e|x|}\frac{\Omega_1(t)}{t}dt.
$$
Clearly, $\Omega_1$ and $\Omega_2$ are nondecreasing as well as $\Omega$
and
$\Omega_2(x)\ge \Omega_1(x)\ge \Omega(x)$. 
We have also $\Omega_2(x)=o(x^{-1+1/\beta})$, $x \to\infty$,
and, consequently, $\int_\mathbb{R} \Omega_2(x) (1+x^2)^{-1}dx<\infty$.
It is shown in \cite{hm2}, Lemma 4.7
that $\widetilde{\Omega}_2$ is a smooth function
and its derivative is given by
$$
\frac{d\widetilde{\Omega}_2(x)}{dx}=-\frac{1}{\pi}
\int\limits_0^\infty \log\left|\frac{1+t}{1-t}\right|
\frac{\Omega(e^2xt)}{|x|t}dt, \qquad x\ne 0.
$$

Since $\Omega(t)=o(t^{-1+1/\beta})$, 
$t\to\infty$, by the dominated convergence theorem, we have
$$
|x|^{2-1/\beta}\frac{d\widetilde{\Omega}_2(x)}{dx}
=-\frac{1}{\pi} \int\limits_0^\infty 
\frac{1}{t^{2-1/\beta}} \log\left|\frac{1+t}{1-t}\right|
\cdot\frac{\Omega(e^2xt)}{(|x|t)^{-1+1/\beta}}dt \to 0
$$                                             
when $|x|\to\infty$. Thus, 
$(\widetilde{\Omega}_2)'(x)=o(|x|^{-2+1/\beta})$, $|x|\to\infty$.
\smallskip

To apply Corollary 2.3 we have to compare the growth of
$\varphi_\beta=\arg B_\beta$ with the growth of $\widetilde{\Omega}_2$.
It follows from the estimates (\ref{17a})--(\ref{17b}) that
$(\widetilde{\Omega}_2)'(x)=o(\varphi'_\beta(x))$ for sufficiently large 
$|x|$. By Corollary 2.3, $w_2=e^{-\Omega_2} \in{\rm Adm}(B_\beta)$
and, consequently, $w\in{\rm Adm}(B_\beta)$. $\bigcirc$
\medskip
\\
{\bf Corollary 5.4.} {\it $\alpha(\beta)\ge \max(1/2, -1+1/\beta)$ 
for $1/2<\beta<1$.}
\smallskip
\\
{\bf Proof.} By Lemma 5.3, $w_\gamma\in{\rm Adm}(B_\beta)$ for $\gamma<-1+1/\beta$.
On the other hand, $\alpha(\beta)\ge\alpha(1)=1/2$. $\bigcirc$
\medskip
\\
{\bf Remark.} Note that the method of Lemma 5.3 
does not allow to prove Theorem 1.1. Indeed,
in the case $\beta=1$ we have $\varphi'_1(x)\asymp |x|^{-1}$, $x<-1$,
whereas for the majorant $w=w_\gamma$ with $\gamma\in (0,1/2)$ 
$$
(\widetilde{\Omega}_2)'(x)\le -C|x|^{\gamma-1},\qquad x<-1,
$$
for some $C>0$. Therefore the function 
$\varphi+ 2\widetilde{\Omega}_2$ is not increasing. Thus, for the case
of sparse zeros the sufficient condition of Corollary 2.3
is far from being necessary. 
\bigskip
\\
{\bf Proof of (\ref{4}) for the case $1/2<\beta<1$.}
By Corollary 5.4, $\alpha(\beta)\ge-1+1/\beta$ for $\beta<1$. Also we have
$\alpha(\beta)\ge1/2$, $\beta<1$, since the function $\alpha(\beta)$
is nonincreasing. Note that $-1+1/\beta=1/2$ for $\beta=2/3$.
\smallskip

Let $1/2<\beta<2/3$. We show that in this case $\alpha_+(\beta)
\le -1+1/\beta$ and, since $\alpha_+(\beta)$ is a nonincreasing
function of $\beta$, the proof of the formula (\ref{4}) will be completed.
\smallskip

Assume that $\gamma>-1+1/\beta$ and there is a nonzero function $f$ 
in ${K_{B_\beta}}$ such that $|f(t)|\le e^{-t^\gamma}$, $t>0$.
Let us show that $f$ is bounded in the domain $\Delta$
defined by \smallskip 
(\ref{dom}).

By (\ref{17b}), $B'_\beta$ is bounded on $(-\infty, 0]$
and we have $f\in H^\infty(\mathbb{C}^+)$ by (\ref{estim1}). 
Applying Lemma 3.4 to the function
$g(z)=f(z)/B_\beta(z)$ in the lower half-plane $\mathbb{C}^-$ 
($g\in H^2(\mathbb{C}^-)$)
we see that 
$$
\log\left|\frac{f(z)}{B_\beta(z)}\right|\le -C |z|^\gamma, 
\qquad -\pi/2\le \arg z <0,
$$
where $\arg z$ stands for the main branch of the argument.
Now it follows from (\ref{l32}) that
$$
|f(z)|\le C_0, \qquad z\in\Delta\cap \{-\pi/2\le \arg z\le 0\}.
$$
By Lemma 3.2, $\log|f(z)|\le C_1+C_2|z|^{-1+1/\beta}$, $z\in\Delta$.  
Applying the Phragm\'{e}n--Lindel\"of principle to the function $f$
in the angle $\{-\pi < \arg z <-\pi/2\}$ we see that $f$
is bounded in $\Delta$. Hence, the function $g=f\circ\eta$, where
$\eta$ is the conformal mapping (\ref{conf}) of $\mathbb{C}^+$ onto $\Delta$,
is bounded in $\mathbb{C}^+$. On the other hand,
$\eta(t)\asymp t^2$, $t\to\infty$, and, consequently,
$\log|g(t)|\le -C_3 t^{2\gamma}$, $t>1$. Note that $2\gamma>1$
since $\gamma>-1+\frac{1}{\beta}$ and $\beta<2/3$. Hence, $g\equiv 0$
and we got a contradiction. $\bigcirc$ 
\medskip
\\
{\bf Remark.} 
It is interesting to compare the formula for $\alpha_+(\beta)$
with the results of A.~Borichev and M. Sodin 
on weighted polynomial approximation on discrete 
subsets of the line (\cite{bs}, Appendix 2).
Let $x_n=n^\beta$, $n\in\mathbb{N}$, and let $w_{\gamma,A}(t)=
\exp(-A|t|^\gamma)$, where $A,\gamma>0$. 
Consider the space
$$
\ell^2(w_{\gamma,A})=\{f:\{x_n\}\to \mathbb{C}: \sum\limits_{n=1}^\infty |f(x_n)|^2
w_{\gamma,A}(x_n)<\infty\}.
$$

The following theorem answers the question about the density of 
the polynomials
in the spaces $\ell^2(w_{\gamma,A})$.
\smallskip
 
{\it If $\beta>2$, then the polynomials are dense in the space
$\ell^2(w_{\gamma,A})$  for $\gamma>1/\beta$ and are not dense for
$\gamma<1/\beta$;
if $\gamma=1/\beta$, then the polynomials are dense if and only if 
$A\ge 2\pi \cot \frac{\pi}{\beta}$. 
If $\beta\le 2$, then the polynomials are dense 
in $\ell^2(w_{\gamma,A})$ if and only if $\gamma \ge 1/2$.}
\smallskip

Thus, for $\beta\ge 2/3$ (but not for $\beta<2/3$) 
the limit exponent $\alpha_+(\beta)$
coincides with the limit $\gamma$ in the theorem of Borichev and Sodin.
Moreover, one can deduce the formula for $\alpha_+(\beta)$
from this theorem
making use of the representations of the form (\ref{8})
with rapidly decreasing $|c_n|$
and an argument analogous to Lemma 3.5. Here we prefer to use
a more direct approach.
\medskip
\\
{\bf Proof of (\ref{5}) for the case $1/2<\beta<1$.}
We use once more the smoothing technique of Lemma 5.3. 
Let $0<\alpha<1$. Consider the functions
\begin{equation}
\label{u1}
U(t)=
\begin{cases}
0, & |t|\le 1, \\
|t|^\alpha -1,  & |t|>1, 
\end{cases}
\end{equation}
and
$$
V(t)=
\begin{cases}
0, & t \le 1, \\
1-t^\alpha,  & t>1. 
\end{cases}
$$
Now let
\begin{equation}
\label{u2}
U_1(x)=\int\limits_{0}^{|x|}\frac{U(t)}{t}dt,
\qquad
U_2(x)=\int\limits_{0}^{|x|}\frac{U_1(t)}{t}dt.
\end{equation}
Analogously, for $x>0$ let 
$$
V_1(x)=\int\limits_{0}^{x}\frac{V(t)}{t}dt,
\qquad
V_2(x)=\int\limits_{0}^{x}\frac{V_1(t)}{t}dt.
$$
For $x<0$ let $V_1(x)=V_2(x)=0$.

We will show that there exist positive constants $K$ and $M$ such that
the majorant
$$
w=\exp(-K U_2 - MV_2)
$$
is in ${\rm Adm}_-(B_\beta)$ for any $\beta\in (1/2,1)$. Clearly,
$K U_2(t) + MV_2(t)\asymp |t|^\alpha$, $t<-1$. Thus 
we obtain the estimate $\alpha_-(\beta) \ge \alpha$. 
Since $\alpha$ is an arbitrary
number from the interval $(0,1)$, we have $\alpha_-(\beta) =1$ 
and our statement will be proved.

By Lemma 4.7 of \cite{hm2},
$$
\frac{d \widetilde{U}_2(x)}{dx}=-\frac{1}{\pi}
\int\limits_0^\infty \log\left|\frac{1+t}{1-t}\right|
\frac{U(xt)}{|x|t}dt, \qquad x\ne 0.
$$
Hence,
\begin{equation}
\label{u3}
\frac{d \widetilde{U}_2(x)}{dx}=-\left(\frac{1}{\pi}
\int\limits_0^\infty \log\left|\frac{1+t}{1-t}\right|
t^{\alpha-1} dt\right) |x|^{\alpha-1}+
O\left(\frac{1}{|x|}\right), \qquad |x|>1.
\end{equation}

Analogously, it is easy to show that
$$
\frac{d \widetilde{V}_2(x)}{dx}=
-\frac{1}{\pi}\int\limits_0^\infty \log\left(\frac{1+t}{t}\right)
\frac{V(|x|t)}{|x|t}dt,  \qquad x< 0,
$$
and
$$
\frac{d \widetilde{V}_2(x)}{dx}=
\frac{1}{\pi}\int\limits_0^\infty \log\left|\frac{1-t}{t}\right|
\frac{V(xt)}{xt}dt,  \qquad x>0.
$$
Hence, for $x<-2$ we have
$$
\frac{d \widetilde{V}_2(x)}{dx}=
\left(\frac{1}{\pi}\int\limits_0^\infty \log\left(\frac{1+t}{t}\right)
t^{\alpha-1} dt\right)|x|^{\alpha-1}
-  
\left(\frac{1}{\pi}\int\limits_0^{1/|x|}\log\left(\frac{1+t}{t}\right)
t^{\alpha-1} dt\right)|x|^{\alpha-1}
$$
$$
-\frac{1}{\pi |x|}\int\limits_{1/|x|}^\infty \log\left(\frac{1+t}{t}\right)
\frac{dt}{t}=
\left(\frac{1}{\pi}\int\limits_0^\infty \log\left(\frac{1+t}{t}\right)
t^{\alpha-1} dt\right)|x|^{\alpha-1}+
O\left( \frac{\log^2 |x|}{|x|}\right). 
$$
Finally, it is easy to see that
$$
\frac{d \widetilde{V}_2(x)}{dx}=
O(x^{\alpha-1}),  \qquad x>2.
$$

Now we take two positive constants $K$ and $M$
such that
$$
K \int\limits_0^\infty \log\left|\frac{1+t}{1-t}\right|
t^{\alpha-1} dt =
M \int\limits_0^\infty \log\left(\frac{1+t}{t}\right)
t^{\alpha-1} dt.
$$
Hence, the function $\Omega=KU_2+MV_2$  satisfies the following 
asymptotic equalities:
$$
\frac{d\widetilde{\Omega}(x)}{dx}=O\left( \frac{\log^2 |x|}{|x|}\right),  
\qquad x<-2;
$$
$$
\frac{d\widetilde{\Omega}(x)}{dx}=O(x^{\alpha-1}),  \qquad x>2.
$$

To apply Corollary 2.3 we should compare the growth of
$\widetilde{\Omega}$ with the growth of the argument $\varphi_\beta$
of the product $B_\beta$. By (\ref{17a}) and (\ref{17b}), 
using that $\beta<1$, we obtain
$$
(\widetilde{\Omega})'(x)=o(\varphi_\beta'(x)),\qquad |x|\to\infty.
$$
Moreover, $\varphi_\beta +2\widetilde{\Omega}$
is a mainly increasing function. Indeed,
$\varphi_\beta +2\widetilde{\Omega}$ is Lipschitz and increasing
on $(-\infty, -R)$ for some $R>0$. To show that
$\varphi_\beta +2\widetilde{\Omega}$ is mainly increasing on 
$(0,\infty)$, put $d_n=(Mn)^{\beta}$
for a sufficiently large $M$.
Hence, $w=e^{-\Omega}$ is admissible, which completes 
the proof of the theorem. $\bigcirc$
\medskip
\\
{\bf Remark 5.5.}
It should be noted that nonzero functions in $K_\beta$, $1/2 < \beta <1$, 
with fast decrease at $-\infty$ should be necessarily unbounded on 
$[0,\infty)$.
For example, any nonzero function $f\in K_{B_\beta}$, $2/3<\beta<1$, 
which is majorized on $(-\infty,0]$ by $w_\alpha$ with 
$\alpha\in (1/2,1)$, is unbounded on $[0,\infty)$. 
\smallskip

Indeed, assume that $|f(t)|\le w_\alpha(t)$, $t<0$, and
$|f(t)|\le 1$, $t>0$. Put $g=f\circ \eta$ where 
$\eta$ is the conformal mapping (\ref{conf}). Then, by Lemmas
3.2, 3.3 and by the arguments analogous to those in the proof 
of Theorem 1.2, we have $|g(t)|\le C_1$, $t>0$,
$\log |g(z)| \le C_2|z|^{-2 + 2/\beta}$, $z\in \mathbb{C}^+$,
and $\log |g(iy)| \le -C_3 y^{2\alpha}$, $y>0$.
Note that $-2 + 2/\beta <1$, since $2/3<\beta<1$,
and $2\alpha>1$. Now $g\equiv 0$ by a variant of the 
Phragm\'{e}n--Lindel\"of principle.
\bigskip

We conclude this section with the formula for the fastest possible
decay of elements of $K_B$ for two-sided sequences with power growth.
Let $B$ be the Blaschke product with the zeros
$$
z_n=
\begin{cases}
n^\beta+i, & n\in \mathbb{Z},\  n>0, \\
-|n|^\gamma+i,  & n\in \mathbb{Z},\  n<0. 
\end{cases}
$$
where $\beta,\gamma>1/2$, and let
$$
\alpha(\beta,\gamma)=\sup\{\alpha: w_\alpha\in{\rm Adm}(B)\}.
$$
\smallskip
\\
{\bf Theorem 5.6.} {\it Let $\beta\ge\gamma >1/2$. If $\beta\le 1$, then 
$\alpha(\beta,\gamma)=1$. If $\beta>1$, then
$$
\alpha(\beta,\gamma)=\max\left(\frac{1}{\beta},\alpha(\gamma)\right).
$$}
{\bf Proof.} The case $\beta\le 1$ is obvious. Let $\beta>1$ and
let $\rho> \max (\frac{1}{\beta},\alpha(\gamma))$. Assume that
$f\in K_B$ and $|f(t)|\le e^{-|t|^\rho}$, $t\in\mathbb{R}$. 
By Lemma 3.4, $|f(z)|\le e^{-C|z|^\rho}$, $z\in\mathbb{C}^+$.
\smallskip

Consider the function $E_\beta$ defined as in (\ref{18})
and put $g=fE_\beta$. It is clear that the function $E_\beta$ 
is of order at most $1/\beta$. Hence, 
$$
|g(z)|\le e^{-C_1|z|^\rho}, \qquad z\in\mathbb{C}^+,
$$
for some $C_1>0$. Thus, $g\in H^2$ and, 
since $g$ is meromorphic in $\mathbb{C}$
and all its poles are in the set $-|n|^{\gamma}-i$, $n< 0$,
we conclude that $g\in K_{B^\#_\gamma}$, where $B^\#_\gamma$
is the Blaschke product with the zeros
$-|n|^\gamma+i$, $n< 0$.  
\smallskip

By our assumption $\rho>\alpha(\gamma)$ and 
$|g(t)|\le e^{-C_1|t|^\rho}$, $t\in\mathbb{R}$. It follows that $g\equiv 0$.
So we see that $w_\rho\notin{\rm Adm}(B)$ when
$\rho> \max (\frac{1}{\beta},\alpha(\gamma))$. Hence,
$\alpha(\beta,\gamma)\le \max (\frac{1}{\beta},\alpha(\gamma))$
\smallskip

Since $K_{B^\#_\gamma}\subset K_B$ it follows that 
$\alpha(\beta,\gamma)\ge \alpha(\gamma)$. We show that 
$\alpha(\beta,\gamma)\ge 1/\beta$, $\beta>1$.
Let $\alpha(\gamma)< 1/\beta$ and $\rho<1/\beta$, and let 
$U$ be the function (\ref{u1}) with $\alpha=\rho$.
Applying the smoothing procedure (\ref{u2}), we get
the majorant $e^{-U_2} \le e w_\rho$, and, by (\ref{u3}),
$(\widetilde{U}_2)'(x)\asymp |x|^{\rho-1}$, $|x|>1$.
\smallskip

Denote by $\varphi$ an increasing branch of the argument of $B$.
Clearly,
$$
\varphi'(x) = \varphi'_\beta(x) +\varphi'_\gamma(-x), \quad x\in\mathbb{R}.
$$
Let us show that the function $\varphi+ 2 \widetilde{U}_2$
is mainly increasing. Let $d_n =n^\beta$, $n\in\mathbb{N}$.
Since $\alpha(\gamma)< 1/\beta$, we have
$1/\gamma<1+1/\beta$. It follows from (\ref{17})--(\ref{17b})
that $\varphi(d_{n+1})-\varphi(d_n)\asymp 1$, whereas 
for $\rho <1/\beta$
$$
\sup\limits_{x_1,\,x_2\in [d_n, d_{n+1}]} 
|\widetilde{U}_2(x_1)-\widetilde{U}_2(x_2)|\le
C((n+1)^{\rho\beta}-n^{\rho\beta})\to 0, \quad n\to\infty.
$$
Thus, we have a sequence $\{d_n\}$ satisfying the conditions
of Theorem 2.2.
On the negative semiaxis $\varphi$ grows even faster than 
for $x>0$ since $\gamma\le\beta$. Hence, for $x<0$ 
the conditions of Theorem 2.2 are also satisfied
(one can take $d_n=-|n|^{\gamma}$, $n<0$) 
and so the function $\varphi+\widetilde{U}_2$ is mainly increasing.
By Theorem 2.2, $e^{-U_2} \in{\rm Adm} (B)$ and, consequently, the majorant
$w_\rho$ is admissible for $K_B$ whenever $\rho<1/\beta$. Therefore 
$\alpha(\beta,\gamma)\ge 1/\beta$, which completes the proof. $\bigcirc$
\smallskip
\\
{\bf Remark.} The case $\beta=\gamma<1$ is considered in more 
detail in \cite{hm2} where certian conditions
sufficient for admissibility are given, in terms
of $\widetilde{\Omega}$. Some admissibility criteria for
two-sided sequences of zeros $z_n$, $n\in\mathbb{Z}$,
are also obtained in \cite{belov} where
the results are stated in terms of $\Omega$ (not $\widetilde{\Omega}$)
and some oscillating $\Omega$ are studied.

\bigskip

\begin{center}
{\bf \large \S 6. Tangential zeros}
\end{center}

In this section we prove Theorem 1.5. Recall that now $B$
is the Blaschke product with zeros $z_n=n+iy_n$, $n\in\mathbb{Z}$,
where $0<y_n\le 1$, the sequence $y_n$ is even and nonincreasing
for $n\ge 0$. 
\smallskip

In what follows we will need the functions $y(t)$, $Y(t)$, 
$t\in\mathbb{R}$, defined in Theorem 1.5. Clearly, the integral
$$
{\cal L}(y)=\int\limits_\mathbb{R}\frac{Y(t)}{1+t^2}dt 
$$
converges if and only if the series (\ref{6}) converges.
\medskip
\\
{\bf Proof of Statement 1 of Theorem 1.5.}
Let majorant $w$ be even, nonincreasing on $\mathbb{R}_+$, and let
${\cal L}(w)<\infty$ (without loss of generality we assume that $w\le 1$). 
We will show that $w\in{\rm Adm}(B)$ under condition (\ref{6}). 
\smallskip

Let us consider the Blaschke product $B^\#$
with the zeros $\zeta_n=n+i(y_n+1)$. It is clear
that $(\arg B^\#)'\asymp 1$ (by the $\arg B^\#$ 
we mean an increasing branch). Hence, by Corollary 2.5,
${\rm Adm}(B^\#)\supset {\rm Adm}(e^{ibz})$ for some $b>0$. 
In particular, the majorant 
\begin{equation}
\label{24c}
w_1(t)=y^2(t)w^{2A}(t+1)(1+t^2)^{-1}
\end{equation}
is in ${\rm Adm}(B^\#)$ for any $A>1$, since ${\cal L}(w_1)<\infty$. 
Thus, there exists a nonzero function $g\in K_{B^\#}$
such that
$$
|g(t)|\le y^2(t)w^{2A}(t+1)(1+t^2)^{-1},\qquad t\in\mathbb{R}.
$$

The sequence $\{\zeta_n\}$ is interpolating. Therefore
$g$ may be represented as
$$
g(z)=\sum\limits_{n\in\mathbb{Z}}\frac{c_n}{z-n+i(y_n+1)}.
$$
Moreover, by Lemma 3.4 and a form of Lemma 3.5, we have 
$$
|c_n|\le C_1 y_n w^A(n+1)(n^2+1)^{-1/2}.
$$
Put $d_n=c_n/\sqrt{y_n}$ and consider the function
$$
f(z)=g(z-i)=\sum\limits_{n\in\mathbb{Z}}\frac{\sqrt{y_n}d_n}{z-n+iy_n}.
$$
Since the sequence $z_n=n+iy_n$ is interpolating and $\{d_n\} \in\ell^2$,
the function $f$ belongs to $K_B$. It remains to verify that
$|f(t)|\le Cw(t)$, $t\in\mathbb{R}$.  Note first that
$$
\left|\frac{c_n}{t-n+iy_n}\right|\le C_1 w^A(n+1), \qquad
t\in \mathbb{R}.
$$
After that, we note that the function
$$
g_n(z)=g(z)-\frac{c_n}{z-n+iy_n}
$$
is analytic in $\Omega_n=\left\{z:n-2/3\le{\rm Re}\, z\le
n+2/3,-1\le{\rm Im}\, z\le 1\right\}$ and,
by (\ref{estim1})--(\ref{estim2}), we have
$|g_n(z)|\le C_2$, $z\in\Omega_n$.
At the same time,   
$$
|g_n(t+i)|\le w^A(t+1)+C_1 w^A(n+1), \qquad t\in[n-2/3,n+2/3].
$$
Therefore, by the theorem on two constants,
$|g_n(t)|\le C_3w(t)$, $t\in[n-1/2,n+1/2]$, for sufficiently
large constant $A$. $\bigcirc$
\bigskip
\\
{\bf Proof of Statement 2 of Theorem 1.5.}
Now we have an additional assumption that 
the function $Y(e^x)$ is a convex function of $x$.
\smallskip

Assume that the majorant $w$ decays faster than any power, that is,
$w(t)=o(|t|^{-N})$, $|t|\to\infty$, for any $N>0$.
Let $f$ be a function in $K_B$ such that $|f|\le w$ on $\mathbb{R}$.
Since the sequence $\{z_n\}$ is interpolating we have
the representation
$$
f(z)=\sum\limits_{n\in\mathbb{Z}}\frac{\sqrt{y_n}c_n}{z-n+iy_n},
$$
where $\{c_n\}\in\ell^2$. Proceeding as in the proof of Lemma 3.5,
we obtain the equalities
\begin{equation}
\label{24a}
\sum\limits_{n\in\mathbb{Z}} \sqrt{y_n}c_n (n-iy_n)^k=0,\qquad k \in \mathbb{Z}_+. 
\end{equation}

Consider the function $F(x)=\int_{\mathbb R}f(t)e^{-itx}dt$,
$x\in\mathbb{R}$. Then 
\begin{equation}
\label{newform}
F(x)=G(x)=-2\pi i \sum\limits_{n\in\mathbb Z}\sqrt{y_n}c_n
\exp[-ix(n-iy_n)], \qquad x>0.
\end{equation}
Since $y_n\le 1$ it follows from the divergence of (\ref{6})
and from the convexity of $Y(e^x)$ that the function
$G$ defined by (\ref{newform}) is $C^\infty$ smooth on $\mathbb{R}$ and
$$
|G^{(k)}(x)|\le C \sup\limits_{n\in\mathbb{N}} 
\sqrt{y_n}\,n^{k}e^k, \qquad x\in\mathbb{R}.
$$
It follows from (\ref{24a}) that
$G^{(k)}(0)=0$, $k\ge 0$. 
Furthermore, 
$$
A_k=\sup\limits_{n\in\mathbb{N}} \sqrt{y_n}\,(en)^{k}\le
\exp\left(\sup\limits_{r>0} [k\log r-
\frac{1}{2}Y(r/e)]\right).
$$
Put $T(r)=\sup_{k\in\mathbb{Z}_+}r^k/A_k$. Then, by a Legendre transform argument
analogous to that in the proof of Lemma 3.6, we get
$$
\log T(r)\ge \frac{1}{2} Y(r/e)-C\log r,\qquad r>1,
$$
for some constant $C$.
Therefore $\int_1^\infty r^{-2} \log T(r)dr=\infty$,
and the classical Denjoy--Carleman quasianalyticity theorem implies
$G\equiv 0$. Since $f\in H^1$ we have 
$$
\int_\mathbb{R} f(t)e^{-itx}dt=0, \qquad x\le 0. 
$$
We conclude that $F\equiv 0$
and, hence, $f\equiv 0$. $\bigcirc$
\bigskip

Results analogous to Theorem 1.5 may be obtained for the
Blaschke products with power growth of zeros
or with one-sided zeros. Let us state
the corresponding result for the case of the Blaschke product
$B_1$ with the zeros $z_n=n+iy_n$, $n\in\mathbb{N}$.
\medskip
\\
{\bf Theorem 6.1.} {\it Let $\{y_n\}_{n\in\mathbb{N}}$
be a positive nonincreasing sequence.
\par
1. If 
\begin{equation}
\label{24}
\sum\limits_{n\in\mathbb{N}} n^{-3/2} \log\frac{1}{y_n}<\infty,
\end{equation}        
then any 
even majorant $w$ nonincreasing on $\mathbb{R}_+$ 
with convergent integral (\ref{3})
is admissible for $K_{B_1}$. 
\smallskip
\par
2. Let $y:[0,\infty) \to (0,\infty)$
be a nonincreasing function such that $y(n)=y_n$, $n\in\mathbb{N}$,
and let $Y=-\log y$.
If the function $Y(e^x)$ is convex on $\mathbb{R}$ 
and the series (\ref{24}) diverges, then any majorant $w$ 
such that $w(t)=o(|t|^{-N})$, $t\to\infty$, 
for every $N>0$, is not admissible for $K_{B_1}$.}
\smallskip
\\
{\bf Proof.} The proof is analogous to the proof of Theorem 1.5.
To prove Statement 1 we define the Blaschke product $B^\#$
as above, that is,
$B^\#$ is the product with the zeros
$n+i(y_n+1)$, $n\in\mathbb{N}$. By Theorem 1.3, each 
even majorant nonincreasing on $\mathbb{R}_+$ 
with convergent integral (\ref{3}) is admissible for 
$K_{B^\#}$. Now we define $w_1$
by formula (\ref{24c}) and complete the proof as above.
\smallskip

To prove Statement 2 we use the same idea as 
in the proofs of Lemmas 3.5 and 3.6. First we note that
if 
$$
f(t)=\sum\limits_{n\in\mathbb{N}}\frac{\sqrt{y_n}c_n}{t-n+iy_n}
$$
is a function from $K_B$ such that for any $N$
we have $|f(t)|=o(t^{-N})$, $t\to \infty$, 
then
$$
\sum\limits_{n\in\mathbb{N}} \sqrt{y_n}c_n (n-iy_n)^k=0,\qquad k \in \mathbb{Z}_+. 
$$
Consider the function
$$
g(t)=\sum\limits_{n\in\mathbb{N}}\frac{\sqrt{y_n}c_n}{t-n+i(y_n+1)}.
$$
Clearly, $g\in K_{B^\#}$.
On the other hand,
$$
(t+i)^kg(t)-\sum\limits_{n\in\mathbb{N}}\frac{(n-iy_n)^k \sqrt{y_n}c_n}
{t-n+i(y_n+1)}=0.
$$
Therefore,
$$
|g(t)|\le C\inf\limits_{k\in\mathbb{Z}_+}\frac{1}{|t+i|^k}
\sup\limits_{n\in\mathbb{N}}[\sqrt{y_n}|n-iy_n|^k].
$$
Repeating the arguments from the proof of Lemma 3.6
we see that
$|g(t)|\le t^{A}e^{-Y(t)}$, $t>1$,
where $A$ is some positive constant. Since the series (\ref{24})
diverges, we have $\int_1^\infty t^{-3/2}Y(t)dt=\infty$ and 
it follows from Theorem 1.3 that $g\equiv 0$.
$\bigcirc$
\medskip
\\
{\bf Remark.} To obtain another proof of Theorem 6.1,
one may consider the function 
$F(x)=\int_{\mathbb R}f(t)e^{-itx}dt$, $x>0$, and show that 
$F$ extends to $\mathbb{C}^+ \cup \mathbb{R}$ with estimates on the derivatives.
To complete the proof, one should apply a slightly modified
version of the quasianalyticity theorem due to B.I. Korenblum
\cite{kor}. On the other hand, making use of the idea of 
the proof of Theorem 6.1 one can give another proof
of Theorem 1.5, Statement 2.
\smallskip

One more proof of this result can be obtained by using a 
result of M.M. Dzhrbashyan \cite[Theorem 24]{mer} on weighted
polynomial approximation on nowhere dense sets dividing the complex plane.

\bigskip
A.D. Baranov: \\
Saint Petersburg State University,\\
Department of Mathematics and Mechanics,\\
28, Universitetski pr., St. Petersburg,\\
198504, Russia
\bigskip
\\
E-mail: {antonbaranov@netscape.net}
\bigskip
\\
A.A. Borichev: \\
Laboratoire d'Analyse et G\'eom\'etrie, \\
Universit\'e Bordeaux 1, \\
351, Cours de la Lib\'eration, \\ 
33405 Talence, France
\bigskip
\\
E-mail: Alexander.Borichev@math.u-bordeaux1.fr
\bigskip
\\
V.P. Havin: \\
Saint Petersburg State University,\\
Department of Mathematics and Mechanics,\\
28, Universitetski pr., St. Petersburg,\\
198504, Russia
\bigskip
\\
E-mail: havin@VH1621.spb.edu


\begin{thebibliography}{25}


\bibitem {bar99} A. D. Baranov, Differentiation in de Branges spaces 
and embedding theorems, {\it Probl. Mat. Anal.} {\bf 19} (1999),
27--68; English transl. in {\it J. Math. Sci. (New-York)} {\bf 101} (2000),
2, 2881--2913.

\bibitem {bar04} A. D. Baranov, Polynomials in the de Branges 
spaces of entire functions, {\it Ark. Mat.} {\bf 44} (2006), to appear.

\bibitem {bh} 
A. D. Baranov, V. P. Havin,  Admissible majorants for model subspaces,
and arguments of inner functions,
{\it Funktsional. Anal. i Prilozhen.} {\bf 40} (2006), to appear.

\bibitem {belov} Yu. S. Belov, Admissibility criteria for model
subspaces with fast growth of the argument of the generating function,
{\it Zap. Nauchn. Semin. POMI}, to appear. 

\bibitem {bm} A. Beurling, P. Malliavin, On Fourier transforms
of measures with compact support, {\it Acta Math.} {\bf 107} (1962), 291--309.

\bibitem {bs} A. Borichev, M. Sodin, The Hamburger moment problem
and weighted polynomial approximation on discrete 
subsets of the real line, {\it J. Anal. Math.} {\bf 76} (1998), 219--264.

\bibitem {br} L. de Branges, {\it Hilbert spaces of entire functions},
Prentice Hall, Englewood Cliffs (NJ), 1968.
 
\bibitem {cr} J. A. Cima, W. T. Ross,  \textit{The backward shift 
on the Hardy space,} Mathematical Surveys
and Monographs, 79, AMS, Providence, RI, 2000.

\bibitem {hj}
V. Havin, B. J\"oricke, 
\textit{The uncertainty principle in harmonic analysis}, 
Springer-Verlag, Berlin, 1994. 

\bibitem {hm1} V. P. Havin,  J. Mashreghi,
Admissible majorants for model subspaces of $H^2$.
Part I: slow winding of the generating inner function,
{\it Can. J. Math.} {\bf 55}, 6 (2003), 1231--1263.

\bibitem {hm2} V. P. Havin,  J. Mashreghi,
Admissible majorants for model subspaces of $H^2$.
Part II: fast winding of the generating inner function,
{\it Can. J. Math.} {\bf 55}, 6 (2003), 1264--1301.

\bibitem {hmn} V. P. Havin,  J. Mashreghi, F. Nazarov,
Beurling--Malliavin multiplier theorem: the 7th proof,
{\it Algebra i Analiz} {\bf 17} (2005), 5, 3--68.

\bibitem {ko1} P. Koosis, 
\textit{The Logarithmic Integral I}, Cambridge Stud. Adv. Math. 
{\bf 12}, 1988.

\bibitem {ko2} P. Koosis, 
\textit{The Logarithmic Integral II}, Cambridge Stud. Adv. Math. 
{\bf 21}, 1992.

\bibitem {ko3} P. Koosis, 
\textit{Le\c cons sur le Th\'eor\`eme de Beurling et Malliavin},
Les Publications CRM, Montr\'eal, 1996.

\bibitem {ko01} P. Koosis, Estimating polynomials and entire 
functions by using their logarithmic sums over complex sequences,  
{\it St. Petersburg Math. J.} {\bf  13} (2002),  5, 757--789.

\bibitem {kor} B. I. Korenblum,
Quasianalytic classes of functions in a circle,
{\it Dokl. Akad. Nauk SSSR} {\bf 164} (1965), 36--39  (Russian). 

\bibitem {lev} B. Ya. Levin, \textit{Distribution of zeros
of entire functions}, GITTL, Moscow, 1956;
English transl.: Amer. Math. Soc., Providence, 1964;
revised edition: Amer. Math. Soc., 1980.

\bibitem {mer} S. N. Mergelyan, Weighted approximations by polynomials,
{\it Uspehi Mat. Nauk} {\bf 11} (1956), 5, 107--152;
English transl.: AMS Translations, Ser. 2, {\bf 10} (1958), 59--106.

\bibitem {nik} N. K. Nikolski, \textit{Treatise on the Shift Operator},
Springer-Verlag, Berlin-Heidelberg, 1986.

\bibitem {nk12} N. K. Nikolski, \textit{Operators, Functions, and Systems: an
Easy Reading. Vol. 1--2}, Math. Surveys Monogr., Vol. 92--93, 
AMS, Providence, RI, 2002.

\end{thebibliography}
\end{document}